\documentclass[11pt]{article}
\usepackage{amssymb,latexsym,amsmath,amsbsy,amsthm,amsxtra,amsgen}
\numberwithin{equation}{section}

\newfont{\msbm}{msbm10 at 11pt}
\newcommand {\R} {\mbox{\msbm R}}

\newcommand {\N} {\mbox{\msbm N}}

\def\be{\begin{equation}}
\def\ee{\end{equation}}
\def\ba{\begin{align}}
\def\ea{\end{align}}
\def\mn{\bigskip\noindent}

\newtheorem{Theo}{Theorem}[section]
\newtheorem{Lemma}[Theo]{Lemma}

\newtheorem{Prop}[Theo]{Proposition}
\newtheorem{Exm}[Theo]{Example}

\begin{document}
\title{A coalescent model for
the effect of advantageous mutations on the genealogy of a population}
\author{by Rick Durrett\thanks{Partially supported by NSF grants from
the probability program
(9877066 and 0202935) and from a joint DMS/NIGMS initiative to support
research in mathematical biology (0201037).}\kern1em and Jason
Schweinsberg\thanks{Supported by an NSF Postdoctoral Fellowship}}
\maketitle

\begin{abstract}
When an advantageous mutation occurs in a population, the favorable
allele may spread to the entire population in a short time, an
event known as a selective sweep.  As a
result, when we sample $n$ individuals from a population and trace their
ancestral lines backwards in time, many lineages may coalesce almost
instantaneously at the time of a selective sweep.  We show that
as the population size goes to infinity, this process converges to a
coalescent process called a coalescent with multiple collisions.
A better approximation for finite populations can be obtained using
a coalescent with simultaneous multiple collisions.  We also show how
these coalescent approximations can be used to get insight into
how beneficial mutations affect the behavior of statistics that
have been used to detect departures from the usual Kingman's coalescent.
\end{abstract}

\section{Introduction}

Our goal in this paper is to describe the coalescent processes that arise
when we consider the genealogy of a population that is affected by
repeated beneficial mutations.  The starting point for this analysis
will be the continuous-time population model introduced by Moran (1958).
In this model, the population size is fixed at $2N$.  Each individual
independently lives for a time that is exponentially distributed with mean $1$
and then is replaced by a new individual.  The parent of the new
individual is chosen at random from the $2N$ individuals, including the
one being replaced.  Note that we can think of the population as
consisting of $2N$ chromosomes of $N$ diploid individuals, so each member
of the population has just one parent.

Suppose we sample $n$ individuals at random from this population at time zero.
To describe the genealogy of the sample, we will define the ancestral process,
which will be a continuous-time Markov process $(\Psi_N(t), t \geq 0)$ whose
state space is the set ${\cal P}_n$ of partitions of $\{1, \dots, n\}$.
The ancestral process describes the coalescence of lineages as we follow 
the ancestral lines of the sampled individuals backwards in time.
More precisely, $\Psi_N(0)$ is the partition of $\{1, \dots, n\}$
into $n$ singletons, and $\Psi_N(t)$ is the partition of $\{1, \dots, n\}$ such
that $i$ and $j$ are in the same block of $\Psi_N(t)$ if and only if the $i$th
and $j$th individuals in the sample have the same ancestor at time $-Nt$.
It is well-known that the process $(\Psi_N(t), t \geq 0)$
is Kingman's coalescent, a coalescent process introduced by
Kingman (1982).  Kingman's coalescent is a ${\cal P}_n$-valued Markov process
that starts from the partition of $\{1, \ldots, n\}$ into singletons.
All transitions involve exactly two blocks of the partition merging together,
and each such transition occurs at rate one.

Within the last decade, progress has been made on describing the genealogy
of populations in models that allow for natural selection. Krone and
Neuhauser (1997) and Neuhauser and Krone (1997) studied a model in which each
individual can be of type $1$ or $2$.  An individual of type $i$ produces
offspring at rate $\lambda_i$, with $\lambda_2 > \lambda_1$ so that type
$2$ is advantageous.  Each new offspring replaces a randomly chosen
individual from the population, and is the same type as its parent with
probability $1 - u_N$ and the opposite type with probability $u_N$.
Under certain assumptions, they show that the genealogy of a sample from the
population can be described using what they call an ancestral selection graph.
Additional work of Donnelly and Kurtz (1999) and Barton, Etheridge, and Sturm
(2004) has incorporated recombination as well as selection into the model.

The ancestral selection graph arises in the limit as
$N \rightarrow \infty$ in the case of weak selection, where the selective
advantage $\lambda_2/\lambda_1 - 1$ and the mutation rates $u_N$ are $O(1/N)$.
Then, as $N \rightarrow \infty$ the fraction of individuals with the
favored allele can be approximated by a diffusion process.  In this paper,
we consider strong selection, where the selective advantage is $O(1)$.
With strong selection, when a beneficial mutation occurs, there is a
positive probability that the beneficial allele will spread to the
entire population, an event known as a selective sweep.

At the end of a selective sweep, the entire population has the favorable
allele, and every member of the population will trace that favorable allele
back to the individual that had the beneficial mutation that caused the
selective sweep.  However, the genealogy becomes more complicated when
we consider recombination.  Diploid individuals usually do not inherit
an identical copy of one of their parent's chromosomes.  Instead, the
inherited chromosome consists of pieces of each of a parent's two
chromosomes.  Since a chromosome is coming
from two places, we need to consider the genealogy not of an entire
chromosome but of a particular site of interest on the chromosome.
When a selective sweep is caused by a beneficial mutation at a site other
than the site of interest, many individuals may trace their gene
at the site of interest back to the individual that had the beneficial
mutation at the beginning of the selective sweep, while others may trace
their gene at the site of interest to a different ancestor because
of recombination between the two sites on the chromosome.  This effect
was first studied by Maynard Smith and Haigh (1974), who called it
the ``hitchhiking effect.''

As we will show, the typical duration of a selective sweep is only
$O(\log N)$.  Therefore, when we speed up time by a factor of $N$ to define
the ancestral process, the selective sweep takes place almost instantaneously.
Consequently, if we sample $n$ individuals some time after a selective
sweep and define the ancestral process as before,
the ancestral process behaves like Kingman's coalescent until
we get back to the time of a selective sweep.  At that time, many
lineages may coalesce because they get traced back to the individual
with the mutation that caused the selective sweep.  This possibility was
observed by Gillespie (2000), who referred to the resulting coalescent
process as the ``pseudohitchhiking model."  We will show that if selective
selective sweeps happen repeatedly throughout the history of a population
at times of a Poisson process, as proposed by Gillespie (2000), then
under suitable assumptions the ancestral processes will converge as $N \rightarrow \infty$
to a coalescent with multiple collisions,
which is a ${\cal P}_n$-valued Markov process in which many blocks
of the partition can merge at once into a single block.  These coalescent
processes were introduced by Pitman (1999) and Sagitov (1999).

While coalescents with multiple collisions are the limiting coalescent processes as $N \rightarrow \infty$, an improved approximation for finite $N$ can be obtained using a coalescent with simultaneous multiple collisions.  Coalescents with simultaneous multiple collisions, which were introduced by Schweinsberg (2000) and M{\"o}hle and Sagitov (2001), are coalescent processes in which many blocks can merge at once into a single block, and many such mergers can occur simultaneously.  They provide a better approximation than coalescents with multiple collisions in this context because, as noted by Barton (1998), Durrett and Schweinsberg (2004a), and Schweinsberg and Durrett (2004),
multiple groups of lineages can coalesce at the time of a selective sweep.

Coalescents with multiple or simultaneous multiple collisions arise as
limits of ancestral processes in populations that occasionally have very
large families because ancestral lines that go back to an individual
with many offspring will coalesce at the same time.  Coalescents
with multiple collisions arise when a single large family is possible in
a given generation, while coalescents with simultaneous multiple collisions arise when
one generation can contain many large families.  For more details,
see Sagitov (1999, 2003),  M{\"o}hle and Sagitov (2001), and Schweinsberg
(2003).  The results in this paper provide a different
biological application of these coalescent processes.

The rest of this paper is organized as follows.  In section 2, we describe our model for how the population evolves when there can be beneficial mutations.  We state our main result, which is that the genealogy of this process converges to a coalescent with multiple collisions.  In section 3, we present the improved approximation involving a coalescent with simultaneous multiple collisions.  The next two sections are devoted to applications of these results.  In section 4, we discuss how multiple mergers affect the number of segregating sites and pairwise differences in a sample of DNA.  These quantities are used in Tajima's $D$-statistic (see Tajima (1989)), which can be used to detect departures from the standard Kingman's coalescent.  In section 5 we discuss how multiple mergers affect the number of mutations that appear on just a single individual in the sample, which is relevant to the test proposed by Fu and Li (1993) for detecting departures from Kingman's coalescent.  Our results suggest that Fu and Li's test should have less power to detect selective sweeps, at least in large samples, than Tajima's $D$-statistic.  Finally, in section 6, we prove the convergence and approximation theorems stated in sections 2 and 3.

\section{Convergence to a coalescent with multiple collisions}

In this section, we give a precise description of our model of a population that experiences beneficial mutations, and we state our main convergence theorem.  We describe what happens following a single beneficial mutation in subsection 2.1, and we consider recurrent beneficial mutations in subsection 2.2.  Then in subsection 2.3, we state the convergence result and give some examples.

\subsection{The effect of a single beneficial mutation}

In this subsection we describe how the population
evolves after one of the $2N$ individuals experiences
a beneficial mutation.  We will denote
the new favorable allele by $B$ and the other allele by $b$.  We assume
the relative fitnesses of the two alleles are $1$ and $1-s$,
so the $B$ alleles will tend to survive longer.  Immediately
after the mutation, one individual has the $B$ allele and $2N-1$ have the
$b$ allele.  Kaplan, Hudson, and Langley (1989) and Stephen, Wiehe, and
Lenz (1992) proposed modeling the fraction of individuals $p(t)$ with the $B$
allele at time $t$ by using the logistic differential equation
$$\frac{dp}{dt} = sp(1-p).$$
This approach has been popular in simulation studies.  However,
Durrett and Schweinsberg (2004a) showed that this approximation is not
very accurate.  Consequently, we will consider instead a modification to
the Moran model that was studied by Durrett and Schweinsberg (2004a) and
Schweinsberg and Durrett (2004).

At one site, each
chromosome has a $B$ or $b$ allele, but we will be interested in the genealogy
at another neutral site at which all alleles have the same fitness.
As in the Moran model, each individual survives for a time that
is exponentially distributed with mean $1$, and then a replacement is
proposed in which the parent of the proposed new individual is chosen
at random from the $2N$ members of the population.  However, to account
for natural selection, whenever a replacement
of a $B$ chromosome with a $b$ chromosome is proposed, the change is
rejected with probability $s$.  Also, to incorporate recombination into
the model, we say that when a new individual is born, it inherits its
alleles at both sites from the same parent with probability $1-r$.
However, with probability $r$, there is recombination between the two
sites, so the new individual inherits its allele at the neutral site
from its parent's other chromosome.  Because we are treating an
individual's two chromosomes as two separate members of the population,
we model this by saying that, with probability $r$, the new individual
inherits the two alleles from two ancestors chosen independently
at random from the population.

Suppose the beneficial mutation appears on
one chromosome at time $0$, and let $X(t)$ be the number of chromosomes
with the favorable allele at time $t$.  Let $\tau = \inf\{t: X(t) \in
\{0, 2N\} \}$ be the time at which either the $B$ or $b$ allele
disappears from the population.  Suppose we take a random sample of
$n$ individuals from the population at time $\tau$.  Let $\Theta$ be
the partition of $\{1, \dots, n\}$ such that $i$ and $j$ are in the same
block of $\Theta$ if and only if the $i$th and $j$th individuals in the
sample have the same ancestor at time zero when we follow the ancestral
lines associated with the neutral site of interest.  The partition $\Theta$ then
describes how the beneficial mutation affects the genealogy of the sample.
We have the following result concerning the distribution of $\Theta$.
Here $Q_{p,n}$, for $p \in [0,1]$, is the distribution of a random partition $\Pi$ obtained as
follows.  First, define a sequence of independent random variables
$(\xi_i)_{i = 1}^n$ such that $P(\xi_i = 1) = p$ and
$P(\xi_i = 0) = 1-p$ for $i = 1, \dots, n$.  Then define $\Pi$ such that one block
of $\Pi$ consists of $\{i \leq n: \xi_i = 1\}$ and the remaining blocks
of $\Pi$ are singletons.

\begin{Prop}
Fix $n \in \N$, and fix $s \in (0,1)$.  Assume there is a constant
$C'$ such that $r \leq C' /(\log N)$ for all $N$.  Let $\alpha =
r \log(2N)/s$, and let $p = e^{-\alpha}$.
\begin{enumerate}
\item There exists a positive constant $C$, depending continuously on $s$
and $\alpha$ but not depending on $N$, such that $|P(\Theta = \pi|X(\tau)
= 2N) - Q_{p,n}(\pi)| \leq C/(\log N)$ for all $\pi \in {\cal P}_n$.

\item Let $\kappa_0$ be the partition of $\{1, \dots, n\}$ into singletons.
There exists a constant $C$, depending continuously on $s$ and $\alpha$ but
not depending on $N$, such that $P(\Theta \neq \kappa_0 \mbox{ and }X(\tau) = 0)
\leq CN^{-1/2}$.
\end{enumerate}
\label{sweepprop1}
\end{Prop}

Note that in this proposition, the selective advantage $s$ is assumed to
be fixed, but the recombination probability $r$ depends on $N$.  Part 1
of the proposition, which is a restatement of Theorem 1.1 of Schweinsberg and
Durrett (2004), implies that as $N \rightarrow \infty$, the distribution of
$\Theta$, conditional on the event that a selective sweep occurs,
converges to $Q_{p,n}$, where $p$ represents the approximate
fraction of lineages that coalesce at the time of the selective sweep.
Part 2 of the proposition, which we prove in Section 6,
shows that lineages typically do not coalesce when the favorable $B$ allele dies out.  The probability that a selective sweep occurs,
and therefore Part 1 of the proposition applies, is
$s/(1 - (1-s)^{2N})$ (see Durrett (2002) or Schweinsberg and Durrett (2004)).

\subsection{A model with recurrent beneficial mutations}

To model a population in which beneficial mutations can occur repeatedly,
we assume that beneficial mutations at different points
on the chromosome occur at times of a Poisson process.
The selective advantage that these mutations provide and the rate of
recombination between the site of interest and the site of the mutation
will be random.  When there is a beneficial mutation in the population,
the population will evolve as described in the previous subsection.
Between these times, the population will follow the standard Moran model.

To be more precise, we will consider the chromosome to be the line
segment $[-L, L]$.  Our goal will be to describe the genealogy of the
site $0$.  For each $N$, the beneficial mutations will be governed by
a Poisson process $K_N$ on $\R \times [-L, L] \times [0,1]$.
If $(t, x, s)$ is a point in $K_N$, then at time $t$, a mutation, which
provides a selective advantage of $s$, will appear at location $x$
on one of the $2N$ chromosomes.  The intensity measure of $K_N$ will be
$\lambda \times \mu_N$, where $\lambda$ denotes Lebesgue 
measure on $\R$ and $\mu_N$ is a finite measure on
$[-L, L] \times [0,1]$ which governs the rates of beneficial mutations.
The recombination probabilities will be determined by a function
$r_N: [-L, L] \rightarrow [0,1]$.  We assume that $r_N(0) = 0$
and $r_N$ is nonincreasing on $[-L, 0]$ and nondecreasing on $[0, L]$.
Beginning at time $t$, the population will evolve according
to the model described in the previous subsection of a population with
a beneficial allele having selective advantage $s$ and recombination
probability $r_N(x)$.  We let $\tau(t)$ denote the first time
that the beneficial mutation that appears at time $t$ either disappears from
the population or is present on all $2N$ chromosomes.

Let ${\cal T}_N = \{t: (t, x, s) \mbox{ is a point in }K_N \mbox{ for some }
x \mbox{ and }s\}$ be the times at which beneficial mutations are
proposed.
Note, however, that we can not define the evolution of the population 
as explained above
if, for some $t_1, t_2 \in {\cal T}_N$, the intervals $[t_1, \tau(t_1)]$
and $[t_2, \tau(t_2)]$ overlap.  There has been some work in the biology
literature on the question of how a selective sweep is affected by another
selective sweep happening at the same time (see, for example, Barton (1995),
Gerrish and Lenski (1998), and Kim and Stephen (2003)).  However, as we
will show, in our model this overlap occurs too infrequently to have any
affect on our results, so we avoid the issue of defining the population
during periods of overlap by allowing a new beneficial
mutation to occur only when there is no other beneficial mutation currently
in the population.  That is, beneficial mutations will occur at the
times in ${\cal T}_N' = \{t \in {\cal T}_N: \tau(u) < t \mbox{ for all }
u \in {\cal T}_N \mbox{ such that }u < t\}$.  Let $${\cal I}_N =
\bigcup_{t \in {\cal T}_N'} [t, \tau(t)].$$  A beneficial mutation
will be present in the population at time $u$ if and only if
$u \in {\cal I}_N$.  For the intervals in ${\cal I}_N$, the evolution of
the population was defined in subsection 2.1.  For the times in $\R \setminus
{\cal I}_N$, we will say that the population evolves according to the
standard Moran model so that the evolution of the population is
well-defined for all of $\R$.

To define the ancestral process $\Psi_N = (\Psi_N(t), t \geq 0)$, we
sample $n$ of the $2N$ individuals at random from the population at time
zero.  We then define $\Psi_N(t)$ to be the partition of $\{1, \dots, n\}$
such that $i$ and $j$ are in the same block of $\Psi_N(t)$ if and only if
the $i$th and $j$th individuals in the sample got
their allele at location $0$ on the chromosome from the same ancestor at
time $-Nt$.  Note that we are again speeding up
time by a factor of $N$ so that, if
there are no beneficial mutations (i.e. if $\mu_N$ is the zero measure),
the ancestral process $\Psi_N = (\Psi_N(t), t \geq 0)$ is Kingman's coalescent.
When we do have beneficial mutations, the ancestral processes will converge
as $N \rightarrow \infty$, under suitable conditions, to a coalescent
with multiple collisions.

\subsection{The main convergence theorem and examples}

Pitman (1999) introduced coalescents with multiple collisions, in which
many blocks of the partition can merge into one.  These coalescent processes
are in one-to-one correspondence with finite measures $\Lambda$ on $[0,1]$,
and the coalescent process associated with a particular measure $\Lambda$
is called the $\Lambda$-coalescent.  We will consider here only
${\cal P}_n$-valued coalescents because they are what we will need to
approximate the genealogy of a sample of size $n$.  However, the
constructions can be extended, using Kolmogorov's Extension Theorem, to yield
coalescent processes that take their values in the set of partitions of $\N = \{1, 2, \dots\}$.

Suppose $(\Pi_n(t), t \geq 0)$
is the ${\cal P}_n$-valued $\Lambda$-coalescent.  Then $\Pi_n(0)$ is the
partition of $\{1, \dots, n\}$ into singletons.  If $\Pi_n(t)$ has $b$ blocks,
then every possible transition involves merging $k$ of the blocks into one,
where $2 \leq k \leq b$.  Denoting the rate of this transition by
$\lambda_{b,k}$, we have
\begin{equation}
\lambda_{b,k} = \int_0^1 x^{k-2} (1-x)^{b-k} \: \Lambda(dx).
\label{cmcmain}
\end{equation}
If $\Lambda = \delta_0$, where $\delta_0$ denotes a unit mass at zero,
then every transition that involves
two blocks merging into one happens at rate one, and no other transitions
are possible.  Thus, the $\delta_0$-coalescent is Kingman's coalescent.

The theorem below states that when we do have beneficial mutations, the
ancestral processes converge as $N \rightarrow \infty$, under suitable
conditions, to a coalescent with multiple collisions.  The multiple mergers
happen at times of selective sweeps.  
Note that the convergence is in the sense of
finite-dimensional distributions.  Convergence in the stronger Skorohod
topology does not hold because, during the short time intervals
when selective sweeps are taking place, $\Psi_N$ may undergo multiple transitions.

\begin{Theo}
Let $\mu$ be a finite measure on $[-L, L] \times [0,1]$, and let
$r:[-L,L] \rightarrow [0,\infty)$ be a bounded continuous function such that
$r(0) = 0$ and $r$ is
nonincreasing on $[-L, 0]$ and nondecreasing on $[0, L]$.  Suppose that, as
$N \rightarrow \infty$, the measures $N \mu_N$ converge weakly to $\mu$ and
the functions $(\log 2N) r_N$ converge uniformly to $r$.  Let $\eta$ be the
measure on $(0,1]$ such that
$$\eta([y, 1]) = \int_{-L}^L \int_0^1 s 1_{\{e^{-r(x)/s}
\geq y\}} \: \mu(dx \times ds)$$ for all $y \in (0,1]$.  Let $\Lambda$
be the measure on $[0,1]$ defined by $\Lambda = \delta_0 + \Lambda_0$,
where $\Lambda_0(dx) = x^2 \eta(dx)$.  Let $\Pi = (\Pi(t), t \geq 0)$ be
the ${\cal P}_n$-valued $\Lambda$-coalescent.  Then, as $N \rightarrow \infty$,
the finite-dimensional distributions of $\Psi_N$ converge to the
finite-dimensional distributions of $\Pi$.
\label{mainth}
\end{Theo}

Note that in Theorem \ref{mainth}, the recombination
probability is $O(1/(\log N))$.  The function $r$ is assumed
to be monotone on $[-L, 0]$ and $[0, L]$ because the greater
the distance between $0$ and the site of the mutation, the
greater the likelihood of recombination between the two
sites.  Also, the rate of beneficial mutations is
$O(1/N)$, so that the multiple mergers caused by
selective sweeps and the ordinary mergers of two lineages at a time are
happening on the same time scale.  If the rate of selective sweeps were
$o(1/N)$, then the multiple mergers would disappear in the limit.  If selective
sweeps occurred on a faster time scale than $O(1/N)$, then the
multiple mergers would dominate for large $N$ and the limiting coalescent
would have no $\delta_0$ component.  Gillespie (2000) considers this possibility and
proposes that it may explain why observed genetic variation does not appear
to be as sensitive to population size as Kingman's coalescent model predicts.
However, in this paper we focus on the case in which both types of
mergers happen on the same time scale.

We now derive the limiting coalescent with multiple collisions in
two natural examples.

\begin{Exm}
{\em Consider the case in which we are concerned only with
mutations at a single site, all of which have the same selective
advantage.  Fix $\alpha > 0$, and let $\mu_N = \alpha N^{-1} \delta_{(z,s)}$
for some $s \in (0,1]$ and $z \in [-L, L]$.  This
means that beneficial mutations that provide selective advantage $s$
appear on the chromosome at site $z$ at times of a Poisson process.
The measures $N \mu_N$ converge to $\mu = \alpha \delta_{(z,s)}$.
Assume that the recombination functions $r_N$ are defined such that
the sequence $(\log 2N)r_N$ converges uniformly to $r$, and let
$\beta = r(z)$.  Then, for all $y \in (0,1]$, we have
$$\eta([y, 1]) = \int_{-L}^{L} \int_0^1 u 1_{\{e^{-r(x)/u} \geq y\}} \:
\mu(dx \times du) = s \alpha  1_{\{e^{-\beta/s} \geq y\}}.$$  Therefore, $\eta$
consists of a mass $s \alpha$ at $p = e^{-\beta/s}$.  It follows from
Theorem \ref{mainth} that the limiting coalescent process is the
$\Lambda$-coalescent, where $\Lambda = \delta_0 + s \alpha p^2 \delta_p$.
Thus, in addition to the mergers involving just two blocks, we have
coalescence events at times of a Poisson process in which we flip
$p$-coins for each lineage and merge the lineages whose coins come up heads.}
\label{exm1}
\end{Exm}

\begin{Exm}
{\em It is also natural to consider the case in which mutations
occur uniformly along the chromosome.  For simplicity, we
will assume that the selective advantage $s$ is fixed.  Let $\lambda$
denote Lebesgue measure on $[-L, L]$.  Suppose $\mu_N = N^{-1}(\alpha \lambda
\times \delta_s)$, so the measures $N \mu_N$ converge to
$\mu = \alpha \lambda \times \delta_s$.  To model recombination
occurring uniformly along the chromosome, we assume that the functions
$(\log 2N)r_N$ converge uniformly to the function $r(x) = \beta|x|$, so the
probability of recombination is proportional to the distance between the
two sites on the chromosome.  For all $y \in (0,1]$, we have
$$\eta([y, 1]) = \alpha s \int_{-L}^L 1_{\{e^{-r(x)/s} \geq y\}} \: dx
= \alpha s \int_{-L}^L 1_{\{e^{-\beta|x|/s} \geq y\}} \: dx.$$  Since
$e^{-\beta|x|/s} \geq y$ if and only if $|x| \leq -(s/\beta)(\log y)$, we have
$$\eta([y,1]) = \min \bigg\{\frac{-2 \alpha s^2 \log y}{\beta}, \: 2 \alpha s L
\bigg\}.$$  Therefore, for $y \geq e^{-\beta L/s}$, we have
$$\frac{d}{dy} \eta([y,1]) = - \frac{2 \alpha s^2}{\beta y}.$$
Let $c = 2 \alpha s^2/\beta$.  It follows that $\eta$ has a density given by
$g_L(y) = c/y$ for $e^{-\beta L/s} \leq y \leq 1$ and $g_L(y) = 0$ otherwise.
By Theorem \ref{mainth}, the finite-dimensional distributions of the
ancestral processes $\Psi_N$ converge to those of the $\Lambda$-coalescent,
where $\Lambda = \delta_0 + \Lambda_0$ and $\Lambda_0$ has density
$h_L(y) = y^2 g_L(y)$.  Note that as $L \rightarrow \infty$, the
density $h_L(y)$ converges to $h(y)$, where $h(y) = c y$
for $y \in [0,1]$ and $h(y) = 0$ otherwise.  We can think of this as the
limiting coalescent for an infinitely long chromosome.}
\label{exm2}
\end{Exm}

\begin{Exm}
{\em Finally, we show that any $\Lambda$-coalescent with a unit mass at zero can
arise as a limit of ancestral processes in this model.  We first show
how to obtain coalescents of the form $\Lambda = \delta_0 + \Lambda_0$,
where $\Lambda_0$ is a finite measure on $[\epsilon, 1]$
and $0 < \epsilon < 1$.  Note that
in Theorem \ref{mainth}, we have $\Lambda_0(dx) = x^2 \eta(dx)$,
so it suffices to show that $\mu$ and $r$ can be chosen to make $\eta$
an arbitrary finite measure on $[\epsilon, 1]$.  Let $G: [\epsilon, 1]
\rightarrow [0, \infty)$ be any nonincreasing left-continuous
function.  We will choose
$\mu$ and $r$ so that $\eta([y, 1]) = G(y)$ for $\epsilon \leq y \leq 1$
and $\eta([0, \epsilon)) = 0$.
Let $L = -\frac{1}{2} \log \epsilon$, and let $\nu$ be the measure on
$[-L, L]$ such that $\nu([-L,0)) = 0$ and, for $\epsilon \leq y \leq 1$,
$\nu([0, -\frac{1}{2} \log y]) = 2G(y)$.  Suppose $r(x) = |x|$ and
$\mu = \nu \times \delta_{1/2}$.  Then, for $\epsilon \leq y \leq 1$,
\begin{align}
\eta([y,1]) &= \int_{-L}^L \int_0^1 s 1_{\{e^{-r(x)/s} \geq y \}} \:
\mu(dx \times ds) \nonumber \\
&= \frac{1}{2} \int_0^L 1_{\{e^{-2x} \geq y\}} \: \nu(dx) =
\frac{1}{2} \nu([0, -(\log y)/2]) = G(y), \nonumber
\end{align}
as claimed.  Thus, we can get the $\Lambda$-coalescent in the limit
if $\Lambda_0((0, \epsilon)) = 0$.
We can obtain an arbitrary $\Lambda$-coalescent by then taking a limit as
$L \rightarrow \infty$ (or $\epsilon \downarrow 0$) as in Example \ref{exm2}.}
\end{Exm}

\section{Approximation by a coalescent with simultaneous multiple collisions}

A key ingredient in the proof of Theorem \ref{mainth}
is part 1 of Proposition \ref{sweepprop1}.  Part 1 of
Proposition \ref{sweepprop1} says that, up to an error of $O(1/(\log N))$,
we can approximate the effect of a selective sweep on the genealogy by
flipping a $p$-coin for each lineage and merging the lineages whose
coins come up heads.  However, Durrett and Schweinsberg (2004a) observed in
simulations that for $N$ between 10,000 and 1,000,000, the approximation in
Proposition \ref{sweepprop1} works poorly, largely because it is possible for
multiple groups of lineages to coalesce at the time of a selective sweep.
By taking this into account, they were able to give a more complicated 
approximation that works much better in simulations and has an error
of only $O(1/(\log N)^2)$.

Before stating this result, we review Kingman's (1978) paintbox construction
of exchangeable random partitions of $\{1, \dots, n\}$.  Let
$$\Delta = \big\{(x_1, x_2, \dots): x_1 \geq x_2 \geq \dots \geq 0,
\sum_{i=1}^{\infty} x_i \leq 1 \big\},$$
and let $G$ be a probability measure on $\Delta$.  We define a
$G$-partition $\Pi$ of $\{1, \dots, n\}$ as follows.
Let $Y = (Y_1, Y_2, \dots)$ be a $\Delta$-valued random variable with
distribution $G$.  Define a sequence $(Z_i)_{i=1}^n$ to be
conditionally i.i.d. given $Y$ such that $P(Z_i = j|Y) = Y_j$ for all
positive integers $j$ and $P(Z_i = 0|Y) = 1 - \sum_{j=1}^{\infty} Y_j$.
Then define $\Pi$ to be the partition such that
distinct integers $i$ and $j$ are in the same
block if and only if $Z_i = Z_j \geq 1$.  We denote the
distribution of a $G$-partition of $\{1, \dots, n\}$ by $Q_{G,n}$.
Note that if $G$ is a unit mass at $(p, 0, 0, \dots)$, then $Q_{G,n} = Q_{p,n}$.

Next, we define a family of distributions $R(\theta, M)$ on $\Delta$ 
by using a stick-breaking construction.
Let $\theta \in [0,1]$, and let $M$ be a positive integer.  Let
$(W_k)_{k=2}^M$ be independent random variables such that $W_k$ has a
Beta($1, k-1$) distribution.  Let $(\zeta_k)_{k=2}^M$
be a sequence of independent random variables such that $P(\zeta_k = 1) =
\theta$ and $P(\zeta_k = 0) = 1 - \theta$ for all $k$.  For
$k = 2, 3, \dots, M$, let $V_k = \zeta_k W_k$.  To perform the
stick breaking, we first break off a fraction $W_M$ of
the unit interval, then break off a fraction $W_{M-1}$ of what
is left over, and so on until we get down to $W_2$.  For
$k = 2, \dots, M$, the length of the $k$th fragment is
${\tilde Y}_k = V_k \prod_{j=k+1}^M (1 - V_j)$,  and the
length of the first fragment is ${\tilde Y}_1 = \prod_{j=2}^M (1 - V_j)$.
Note that $\sum_{k=1}^M {\tilde Y}_k = 1$.  Let
$Y = (Y_1, Y_2, \dots, Y_M, 0, 0, \dots ) \in \Delta$ be the sequence
obtained by ranking the interval lengths ${\tilde Y}_1, \dots, {\tilde Y}_M$ in
decreasing order and then appending an infinite sequence of zeros.
Finally, let $R(\theta, M)$ be the distribution of $Y$.

These distributions $R({\theta, M})$ were studied in Durrett and Schweinsberg
(2004b), who used them to approximate the distribution of family sizes
in a Yule process with infinitely many types.  They arise in the
proposition below because, after a beneficial mutation, the number of
lineages with the $B$ allele that do not eventually die out can be approximated
by a Yule process.  The result below is Theorem 1.2 of Schweinsberg and Durrett (2004).  

\begin{Prop}
Fix $n \in \N$, and fix $s \in (0,1)$.  Assume there is a constant
$C'$ such that $r \leq C'/(\log N)$ for all $N$.  Let $\alpha =
r \log(2N)/s$, and let $p = e^{-\alpha}$.  Then
there exists a positive constant $C$, depending continuously on $s$
and $\alpha$ but not depending on $N$, such that $$|P(\Theta = \pi|X(\tau) 
= 2N) - Q_{R(r/s, \lfloor 2Ns \rfloor), n}(\pi)| \leq C/(\log N)^2$$
for all $\pi \in {\cal P}_n$, where $\lfloor m \rfloor$
denotes the greatest integer less than or equal to $m$.
\label{sweepprop2}
\end{Prop}

Because the improved approximation allows many groups of lineages to coalesce at the time of a selective sweep, this result suggests that, for finite $N$, a coalescent with simultaneous multiple collisions should provide a better approximation of the ancestral process than a coalescent with multiple collisions.  Coalescents with
simultaneous multiple collisions, which were studied by M{\"o}hle and
Sagitov (2001), Schweinsberg (2000), and Bertoin and Le Gall (2003),
have the property that many blocks can merge at once into
a single block, and many such mergers can occur simultaneously.
Coalescents with simultaneous multiple collisions are in one-to-one
correspondence with finite measures $\Xi$ on $\Delta$.  

Suppose $\pi$ is a partition of $\{1, \dots, n\}$ whose blocks are
$B_1, \dots, B_m$, and suppose $\pi'$ is a partition of $\{1, \dots, n'\}$
with $n' \geq m$ whose blocks are $B_1', \dots, B_k'$.
Following Bertoin and Le Gall (2003), define the coagulation of $\pi$ by $\pi'$
to be the partition whose blocks are given by $\bigcup_{j \in B_i'} B_j$
for $i = 1, \dots, k$.  Suppose
$(\Pi_n(t), t \geq 0)$ is the ${\cal P}_n$-valued $\Xi$-coalescent.
If there are $b$ blocks at time $t-$ and a merger occurs at time $t$,
then there exists a unique partition $\pi \in {\cal P}_b$ such that
$\Pi_n(t)$ is the coagulation of $\Pi_n(t-)$ by $\pi$.  If $\pi$ has
$r+s$ blocks, $s$ of which are singletons and the other $r$ of which have
sizes $k_1, \dots, k_r \geq 2$, where $b = k_1 + \dots + k_r + s$, then
the rate of this transition is
\begin{equation}
\lambda_{b;k_1, \dots, k_r;s} = \int_{\Delta}
Q_{\delta_x, b}(\pi) \bigg( \sum_{j=1}^{\infty} x_j^2 \bigg)^{-1} \: \Xi_0(dx)
+ a1_{\{r = 1, k_1 = 2\}},
\label{csmcmain}
\end{equation}
where $\delta_x$ denotes a unit mass at $x = (x_1, x_2, \dots) \in \Delta$ and
$\Xi$ has been written as $a \delta_{(0, 0, \dots)} + \Xi_0$ with
$\Xi_0(\{(0, 0, \dots )\}) = 0$.
Coalescents with multiple collisions are a special case in which $\Xi$ is
concentrated on points in which only the first coordinate is nonzero.

Coalescents with multiple and simultaneous multiple collisions can be
constructed from Poisson point processes (see Pitman (1999) and
Schweinsberg (2000)).  Consider a Poisson process on $(0, \infty) \times {\cal P}_n$ whose
intensity measure is the product of Lebesgue measure on $(0, \infty)$
and a measure $L$ on ${\cal P}_n$ defined as follows.  Let
$S \subset {\cal P}_n$ be the set of all partitions consisting of one block
of size $2$ and $n-2$ singletons.  If $\pi \in {\cal P}_n$, let $L(\pi) = 0$
if $\pi$ is the partition consisting of $n$ singletons.  Otherwise, let
\begin{equation}
L(\pi) = \int_{\Delta} Q_{\delta_x, n}(\pi) \bigg( \sum_{j=1}^{\infty} x_j^2
\bigg)^{-1} \: \Xi_0(dx) + a 1_{\{\pi \in S\}}.
\label{csmcpois}
\end{equation}
Since $L$ is a finite
measure, it is easy to define $\Pi_n = (\Pi_n(t), t \geq 0)$ such that
$\Pi_n(0)$ is the partition consisting of $n$ singletons and, at the times
of points $(t,\pi)$ of the Poisson point process, the partition
$\Pi_n(t)$ is the coagulation of $\Pi_n(t-)$ by $\pi$, and these are
the only jump times of $\Pi_n$.  This coalescent
process is the ${\cal P}_n$-valued $\Xi$-coalescent.  The construction of
the $\Lambda$-coalescent is the same, except that if $\pi$ has at least
one block that is not a singleton, we define
\begin{equation}
L(\pi) = \int_0^1
Q_{p,n}(\pi) p^{-2} \: \Lambda_0(dp) + a 1_{\{\pi \in S\}},
\label{cmcpois}
\end{equation}
where $\Lambda = \delta_0 + \Lambda_0$ and $\Lambda_0(\{0\}) = 0$.

Under some additional assumptions, most significantly
restricting the selective advantage resulting from each beneficial mutation
to be at least $\epsilon > 0$, we are able to obtain bounds on the
difference between the finite-dimensional distributions of $\Psi_N$ and
the finite-dimensional distributions of the approximating coalescent 
process.  Proposition \ref{coalprop} below shows that indeed the
coalescent with simultaneous multiple collisions gives a more accurate
approximation. 

\begin{Prop}
Let $\mu$ be a finite measure on $[-L, L] \times [\epsilon, 1]$,
where $\epsilon > 0$, and let
$r:[-L,L] \rightarrow [0,1]$ be a function such that $r(0) = 0$ and
$r$ is nonincreasing on $[-L, 0]$ and nondecreasing on $[0, L]$.
Suppose that, for all $N$, we have $\mu_N = N^{-1} \mu$.  Also, assume
that $r_N(x) = r(x)/\log(2N)$ for all $N$ and $x$.  Fix times
$0 < u_1 < \dots < u_m$, and let $\pi_1, \dots, \pi_m \in {\cal P}_n$.
\begin{enumerate}
\item Define $\eta$ and $\Lambda$ as in Theorem \ref{mainth}. 
Let $\Pi = (\Pi(t), t \geq 0)$ be the ${\cal P}_n$-valued $\Lambda$-coalescent.
Then there exists a constant $C$
such that $$|P(\Psi_N(u_i) = \pi_i \mbox{ for }i = 1, \dots, m) -
P(\Pi(u_i) = \pi_i \mbox{ for }i = 1, \dots, m)| \leq \frac{C}{\log N}.$$

\item Let $G_N$ be the measure on $\Delta$ such that for all measurable
subsets $A \subset \Delta$, we have
$$G_N(A) = \int_{-L}^L \int_0^1 s R(r_N(x)/s, \lfloor 2Ns \rfloor)(A)
\: \mu(dx \times ds).$$  Let $\Xi_N$ be the measure on $\Delta$ given by
$\Xi_N = \delta_{(0, 0, \dots)} + \Xi_{N,0}$, where
$\Xi_{N,0}$ is defined by $\Xi_{N,0}(dx) =
(\sum_{j=1}^{\infty} x_j^2) G_N(dx)$.  Let $\Upsilon_N =
(\Upsilon_N(t), t \geq 0)$ be the ${\cal P}_n$-valued $\Xi_N$-coalescent.
Then there exists a constant $C$
such that $$|P(\Psi_N(u_i) = \pi_i \mbox{ for }i = 1, \dots, m) -
P(\Upsilon_N(u_i) = \pi_i \mbox{ for }i = 1, \dots, m)|
\leq \frac{C}{(\log N)^2}.$$
\end{enumerate}
\label{coalprop}
\end{Prop}

\section{Segregating sites and pairwise differences}

One motivation for modeling a population that experiences recurrent
selective sweeps by coalescents with multiple
or simultaneous multiple collisions is that these coalescent models can
provide insight into tests used to detect selective sweeps.
In view of part 2 of Proposition \ref{coalprop} and the simulation
results in Durrett and Schweinsberg (2004a), there should be little loss
of accuracy in studying the behavior of these tests under the assumption
that the genealogy of a sample follows a coalescent with
simultaneous multiple collisions.
One commonly used test is based on Tajima's $D$-statistic 
(see Tajima (1989)).  Given a sample of $n$ strands of DNA from the same
region on a chromosome, let $\Delta_{ij}$ be the number of sites at which
the $i$th and $j$th segments differ, and let
$\Delta_n = \binom{n}{2}^{-1} \sum_{i \neq j}
\Delta_{ij}$ be the average number of pairwise differences over the
$\binom{n}{2}$ possible pairs.  Let $S_n$ be the number of segregating
sites in the sample, that is, the number of sites at which at least one pair
of segments differs.  Tajima's $D$-statistic compares
the statistics $\Delta_n$ and $S_n$.

Suppose the ancestral history of a sample of $N$ individuals is given
by a coalescent with multiple or simultaneous multiple collisions.
Let $\lambda_b$ be the total rate of all mergers when the coalescent
has $b$ blocks.  Assume that, on the time scale of the coalescent
process, mutations happen at rate $\theta/2$.  Any mutation on the
$i$th or $j$th lineage before these lineages coalesce will cause the
$i$th and $j$th segments to differ at some site.  Since the expected
time for these lineages to coalesce is $\lambda_2^{-1}$, we have
$E[\Delta_{ij}] = \theta \lambda_2^{-1}$.  Therefore
\begin{equation}
E[\Delta_n] = \theta \lambda_2^{-1}.
\label{pairdiff}
\end{equation}
Note that $\lambda_2 = \Lambda([0,1])$ for coalescents with multiple
collisions and $\lambda_2 = \Xi(\Delta)$ for coalescents with simultaneous
multiple collisions.

To calculate the expected number of segregating sites, we note that any
mutation in the ancestral tree before all $n$ lineages have coalesced
into one adds to the number of
segregating sites.  If, at some time, the coalescent has exactly $b$
blocks, the expected time that the coalescent has $b$ blocks is
$\lambda_b^{-1}$.  Let $G_n(b)$ be the probability that
the coalescent, starting with $n$ blocks, will have exactly $b$ blocks
at some time.  Then
\begin{equation}
E[S_n] = \frac{\theta}{2} \sum_{b=2}^n b \lambda_b^{-1} G_n(b).
\label{segsites}
\end{equation}
Although we do not have a closed-form expression for $G_n(b)$, 
these quantities can be calculated recursively because
(\ref{cmcmain}) and (\ref{csmcmain}) allow us to express
$G_n(b)$ in terms of $G_k(b)$ for $k < n$.  As a result, it would
not be difficult to evaluate the expression in (\ref{segsites}) numerically.

Suppose the ancestral process is given by Kingman's coalescent, which
would be the case if there were no selective sweeps.  Then
$\lambda_b = \binom{b}{2}$ for all $b \geq 2$.  Also, the number of
blocks never decreases by more than one at a time, so
$G_n(b) = 1$ whenever $2 \leq b \leq n$.  It follows that
$E[\Delta_n] = \theta$ and
\begin{equation}
E[S_n] = \frac{\theta}{2} \sum_{b=2}^n b \binom{b}{2}^{-1}
= \theta \sum_{b=2}^n \frac{1}{b-1} = \theta h_{n-1},
\label{kingseg}
\end{equation}
where $h_{n-1} = \sum_{i=1}^{n-1} (1/i)$.  Thus,
$E[\Delta_n - S_n/h_{n-1}] = 0$.  This observation is the basis for
Tajima's $D$-statistic, which is given by
\begin{equation}
D = \frac{\Delta_n - S_n/h_{n-1}}{\sqrt{a_nS_n + b_nS_n(S_n-1)}},
\label{tajeq}
\end{equation}
where $a_n$ and $b_n$ are somewhat complicated constants that are
chosen to make the variance of $D$ approximately one when the ancestral
tree is given by Kingman's coalescent.  See section 4.1 of Durrett (2002) for details.

After a selective sweep, the new mutants will tend to have low frequency.
As a result, a recent selective sweep should decrease $\Delta_n$ more
than $S_n$, causing the numerator of Tajima's $D$-statistic to be negative.
Braverman et. al. (1995) found in simulations that Tajima's $D$-statistic
indeed tends to be negative after a selective sweep.  Simonsen, Churchill,
and Aquadro (1995) studied this question further and argued that unless the
selective sweep was recent, Tajima's $D$-statistic had relatively little
power to detect selective sweeps.  
See also Przeworski (2002), who discusses the power of Tajima's $D$-statistic
to detect selective sweeps.
Our coalescent approximation allows us to obtain the following result
regarding the expected number of segregating sites when the population 
experiences recurrent selective sweeps.

\begin{Prop}
Consider a $\Lambda$-coalescent in which $\Lambda = \delta_0 + \Lambda_0$,
where $\Lambda_0(\{0\}) = 0$, or a $\Xi$-coalescent in which
$\Xi = \delta_{(0, 0, \dots)} + \Xi_0$ and $\Xi_0(\{(0, 0, \dots)\}) = 0$.
Let $\alpha_b = \lambda_b - \binom{b}{2}$.
Suppose
\begin{equation}
\sum_{b=2}^{\infty} \frac{\alpha_b \log b}{b^2} < \infty.
\label{segcond}
\end{equation}
Then, there exists a constant $\rho \geq 0$ such that
\begin{equation}
\lim_{n \rightarrow \infty} E[S_n] - \theta h_{n-1} = -\rho.
\label{segeq}
\end{equation}
Furthermore, defining $G_{\infty}(b) = \lim_{n \rightarrow \infty} G_n(b)$,
we have
\begin{equation}
\rho = \frac{\theta}{2}
\sum_{b=2}^{\infty} b \bigg( \binom{b}{2}^{-1} -
\lambda_b^{-1} \bigg) + \frac{\theta}{2} \sum_{b=2}^{\infty} b \lambda_b^{-1}
(1 - G_{\infty}(b)).
\label{rhodef}
\end{equation}
\label{segprop}
\end{Prop}

The condition (\ref{segcond}) prevents
$\Lambda_0$ or $\Xi_0$ from having too much mass near zero.  Note
that (\ref{pairdiff}) implies that $E[\Delta_n]$ decreases
by a constant as a result of the beneficial mutations, while
Proposition \ref{segprop} implies that when (\ref{segcond}) holds, $E[S_n/h_{n-1}]$
decreases by approximately $\rho/h_{n-1}$, which is $O(1/(\log n))$.
Therefore, Proposition \ref{segprop} shows that for
sufficiently large samples we do expect Tajima's
$D$-statistic to be negative when the population is affected by
recurrent selective sweeps.  Before proving this proposition, we consider some examples.

\begin{Exm}
{\em Suppose, as in Example \ref{exm1}, we have a $\Lambda$-coalescent
in which $\Lambda = \delta_0 + s \alpha p^{-2} \delta_p$.  
Since $p$-mergers occur at rate $s \alpha$, we have
$\lambda_b \leq \binom{b}{2} + s \alpha$ and thus $\alpha_b \leq s \alpha$
for all $b$.  Condition (\ref{segcond}) follows immediately.

Suppose instead we have the $\Lambda$-coalescent of Example \ref{exm2},
where $\Lambda = \delta_0 + \Lambda_0$ and $\Lambda_0(dx) = cx \: dx$.
Note that $\alpha_b$ is the same as the total merger rate of the
$\Lambda_0$-coalescent when there are $b$ blocks.  Using the fact that
if $Z \sim \mbox{Binomial}(b,x)$ then $P(Z \geq 2) =
1 - (1-x)^b - bx(1-x)^{b-1}$, we have
\begin{align}
\alpha_b &= \int_0^1 (1 - (1-x)^b - bx(1-x)^{b-1}) x^{-2} \: \Lambda_0(dx)
\nonumber \\
&= c \int_0^1 (1 - (1-x)^b - bx(1-x)^{b-1}) x^{-1} \: dx
\leq c \int_0^1 (1 - (1-x)^b) x^{-1} \: dx \nonumber \\
&= c \int_0^{1/b}  (1 - (1-x)^b) x^{-1} \: dx +
c \int_{1/b}^1 (1 - (1-x)^b) x^{-1} \: dx \nonumber \\
&\leq c \int_0^{1/b} b \: dx + c \int_{1/b}^1 x^{-1} \: dx =
c(1 + \log b),
\label{alphab}
\end{align}
which implies (\ref{segcond}).
}
\end{Exm}

\begin{Exm}
{\em Although (\ref{segcond}) holds in the natural cases given in
Examples \ref{exm1} and \ref{exm2}, we show here that it does not hold
for all coalescents.  Suppose $\Lambda = \delta_0 + \Lambda_0$,
where $\Lambda_0$ is the uniform distribution on $(0,1)$.
Note that there exists a constant $C > 0$ such that if
$Z \sim \mbox{Binomial}(b,x)$ with $x \geq 1/b$ and $b \geq 2$, then
$P(Z \geq 2) \geq C$.  Therefore,
\begin{align}
\alpha_b &= \int_0^1 (1 - (1-x)^b - bx(1-x)^{b-1}) x^{-2} \: dx \nonumber \\
&\geq \int_{1/b}^1 (1 - (1-x)^b - bx(1-x)^{b-1}) x^{-2} \: dx \nonumber \\
&\geq C \int_{1/b}^1 x^{-2} \: dx = C(b-1), \nonumber
\end{align}
so (\ref{segcond}) does not hold in this case.
}
\end{Exm}

\begin{proof}[Proof of Proposition \ref{segprop}]
When the coalescent has $n+1$ blocks, the probability
that the next coalescence event will take the coalescent down to
fewer than $n$ blocks is at most
$[\lambda_{n+1} - \binom{n+1}{2}]/\lambda_{n+1}$.  Therefore,
if $2 \leq b \leq n$, then
\begin{equation}
|G_{n+1}(b) - G_n(b)| \leq
\frac{\lambda_{n+1} - \binom{n+1}{2}}{\lambda_{n+1}}
= \frac{\alpha_{n+1}}{\lambda_{n+1}} \leq \frac{2 \alpha_{n+1}}{n(n+1)}.
\label{Gneq}
\end{equation}
Therefore, when (\ref{segcond}) holds,
the sequence $(G_n(b))_{n=b}^{\infty}$ is Cauchy and thus
has a limit $G_{\infty}(b)$.

It follows from (\ref{segsites}) and (\ref{kingseg}) that
\begin{align}
E[S_n] - \theta h_{n-1} &= \frac{\theta}{2} \sum_{b=2}^n
b \lambda_b^{-1} G_n(b) - \frac{\theta}{2} \sum_{b=2}^n b \binom{b}{2}^{-1}
\nonumber \\
&= \frac{\theta}{2} \sum_{b=2}^n b \bigg( \lambda_b^{-1} - \binom{b}{2}^{-1}
\bigg) - \frac{\theta}{2} \sum_{b=2}^n b \lambda_b^{-1} (1 - G_{\infty}(b))
\nonumber \\
&\hspace{.5in} + \frac{\theta}{2}
\sum_{b=2}^n b \lambda_b^{-1} (G_n(b) - G_{\infty}(b)). 
\label{segmain}
\end{align}
To prove Proposition \ref{segprop}, we need to take the limit as
$n \rightarrow \infty$ of the three terms on the right-hand side 
of (\ref{segmain}).

For the first term, we note that
$$\binom{b}{2}^{-1} - \lambda_b^{-1} =
\frac{\lambda_b - \binom{b}{2}}{\binom{b}{2} \lambda_b} \leq
\alpha_b \binom{b}{2}^{-2} = \frac{4 \alpha_b}{b^2(b-1)^2}.$$
Therefore, when (\ref{segcond}) holds, we have a summable series and
\begin{equation}
\lim_{n \rightarrow \infty} \frac{\theta}{2}  \sum_{b=2}^n b
\bigg( \lambda_b^{-1} - \binom{b}{2}^{-1} \bigg) =
- \frac{\theta}{2} \sum_{b=2}^{\infty} b
\bigg( \binom{b}{2}^{-1} - \lambda_b^{-1} \bigg).
\label{term1}
\end{equation}
For the second term, note that (\ref{Gneq}) and the fact that
$G_b(b) = 1$ imply
\begin{align}
\sum_{b=2}^{\infty} b \lambda_b^{-1} (1 - G_{\infty}(b))
&\leq \sum_{b=2}^{\infty}
\frac{2}{b-1} \bigg( \sum_{m=b}^{\infty} \frac{2 \alpha_{m+1}}{m(m+1)} \bigg)
\nonumber \\
&= \sum_{m=2}^{\infty} \frac{2 \alpha_{m+1}}{m(m+1)}  \sum_{b=2}^m
\frac{2}{b-1} \leq \sum_{m=2}^{\infty}
\frac{4 \alpha_{m+1} (1 + \log (m-1))}{m(m+1)},
\nonumber
\end{align}
which is finite by (\ref{segcond}).  Therefore,
\begin{equation}
\lim_{n \rightarrow \infty} \frac{\theta}{2} \sum_{b=2}^n
b \lambda_b^{-1} (1 - G_{\infty}(b)) = \frac{\theta}{2}
\sum_{b=2}^{\infty} b \lambda_b^{-1} (1 - G_{\infty}(b)).
\label{term2}
\end{equation}
Finally, for the third term,
\begin{align}
\limsup_{n \rightarrow \infty} \sum_{b=2}^n b \lambda_b^{-1}
|G_n(b) - G_{\infty}(b)| &\leq \limsup_{n \rightarrow \infty}
\sum_{b=2}^n \frac{2}{b-1} \bigg( \sum_{m=n}^{\infty}
\frac{2 \alpha_{m+1}}{m(m+1)} \bigg) \nonumber \\
&\leq \limsup_{n \rightarrow \infty}
\frac{1}{\log n} \sum_{b=2}^n \frac{2}{b-1} \bigg(
\sum_{m=n}^{\infty} \frac{2 \alpha_{m+1} \log m}{m(m+1)} \bigg) \nonumber \\
&\leq \limsup_{n \rightarrow \infty}
\frac{2(1 + \log (n-1))}{\log n} \sum_{m=n}^{\infty}
\frac{2 \alpha_{m+1} \log m}{m(m+1)} = 0
\label{term3}
\end{align}
by (\ref{segcond}).  The proposition follows from
(\ref{segmain}), (\ref{term1}), (\ref{term2}), and (\ref{term3}).
\end{proof}

\section{The number of singletons}

Fu and Li (1993) proposed another test to detect departures from Kingman's coalescent.  They considered the ancestral tree in which the leaves are the $n$ individuals in the sample.  They defined the branches connecting a leaf to an internal node to be external branches and the other branches to be internal branches.  Let $\eta_e$ denote the number of mutations on external branches, and let $\eta_i$ be the number of mutations on internal branches.  Every mutation produces a segregating site, so $\eta_e + \eta_i = S_n$.  If a mutation occurs on an external branch, the mutant gene appears on just one of the $n$ individuals in the sample, while if a mutation occurs on an internal branch, the mutant gene appears on between $2$ and $n-1$ of the individuals in the sample.  Therefore, to determine $\eta_e$, we simply count the number of mutations that appear on just one of the sampled chromosomes.  Note that unless an outgroup is available,
it will not be possible to distinguish between a mutation that appears on one of the sampled chromosomes and a mutation that appears on $n-1$ of the sampled chromosomes.  Fu and Li (1993) proposed a modification of their test for when there is no outgroup, but for the analysis in this section, we assume that we have an outgroup that enables us to make this distinction.  

Let $J_n$ be the sum of the lengths of the external branches.  In terms of the associated coalescent process, $J_n$ is the sum, over $i$ between $1$ and $n$, of the amount of time that the integer $i$ is in a singleton block.  Let $I_n$ be the sum of the lengths of the internal branches.  Assuming, as before, that mutations occur at rate $\theta/2$ on the time scale of the coalescent process, we have $E[\eta_e|J_n] = (\theta/2)J_n$ and $E[\eta_i|I_n] = (\theta/2)I_n$.

Fu and Li's $D$-statistic is based on comparing $\eta_i$ with $(h_{n-1} - 1) \eta_e$.  Note that $\eta_i - (h_{n-1} - 1) \eta_e = S_n - h_{n-1} \eta_e$.  To see that this has mean zero when the ancestral tree is given by Kingman's coalescent, we follow the explanation on p. 163 of Durrett (2002).  In the case of Kingman's coalescent, (\ref{kingseg}) gives $E[S_n] = \theta h_{n-1}$.  Therefore,
$E[S_n - h_{n-1} \eta_e] = \theta h_{n-1} - \theta h_{n-1} E[J_n]/2$, so it remains to show that
$E[J_n] = 2$.  Let $K_n$ be the amount of time that the integer $1$ is in a singleton block of the partition, so $E[J_n] = n E[K_n]$.  Let $T_n$ be the amount of time before the first coalescence event, and note that $E[T_n] = 2/[n(n-1)]$.  The probability that $1$ coalesces with another integer at time
$T_n$ is $2/n$, and this event is independent of $T_n$.  If $1$ does not coalesce at this time, then the expected additional time that $1$ is a singleton is $E[K_{n-1}]$.  Therefore, we get the recursion $$E[K_n] = \frac{2}{n} E[T_n] + \frac{n-2}{n} E[T_n + K_{n-1}] = \frac{2}{n(n-1)} + \frac{n-2}{n} E[K_{n-1}].$$  Note that $E[K_2] = 1$, and then it is easy to show by induction that $E[K_n] = 2/n$ for all $n$, and so $E[J_n] = 2$ for all $n$,
as claimed.

We can write Fu and Li's $D$-statistic as
\begin{equation}
D = \frac{S_n - h_{n-1} \eta_e}{\sqrt{c_n S_n + d_nS_n^2}},
\label{fulieq}
\end{equation}
where, as in (\ref{tajeq}), $c_n$ and $d_n$ are constants chosen to make the variance of the statistic approximately one when the genealogy is given by Kingman's coalescent.  Details of the variance computation are given in section 4.2 in Durrett (2002), where an error of Fu and Li (1993) is corrected.

When multiple mergers cause many lineages to coalesce at once, one expects $I_n$ to be reduced more than $J_n$ because there is still an external branch associated with each leaf, but there are fewer internal branches because of multiple mergers.  This would cause Fu and Li's $D$-statistic to be negative.  The next proposition shows that this is indeed the case.

\begin{Prop}
Let $(\Pi_n(t), t \geq 0)$ be a ${\cal P}_n$-valued $\Lambda$-coalescent in which
$\Lambda = \delta_0 + \Lambda_0$, where $\Lambda_0(\{0\}) = 0$, or a ${\cal P}_n$-valued
$\Xi$-coalescent in which $\Xi = \delta_{(0, 0, \dots)} + \Xi_0$ and $\Xi_0(\{(0, 0, \dots)\}) = 0$.
Let $\alpha_b = \lambda_b - \binom{b}{2}$, and suppose (\ref{segcond}) holds.  Then
\begin{equation}
\lim_{n \rightarrow \infty} E[S_n - h_{n-1} \eta_e] = - \rho,
\label{rho2}
\end{equation}
where $\rho$ is the constant defined in (\ref{rhodef}).
\label{fuliprop}
\end{Prop}

The key to the proof of this proposition is the following lemma.

\begin{Lemma}
Under the assumptions of Proposition \ref{fuliprop}, there is a positive constant $C$ such that
\begin{equation}
0 \leq E[2 - J_n] \leq \frac{C}{n} \sum_{b=2}^n \frac{\alpha_b}{b}
\label{Lneq}
\end{equation} 
for all $n \geq 2$.
\label{singlem}
\end{Lemma}

The first inequality in (\ref{Lneq}), which does not require condition (\ref{segcond}), shows that the expected sum of the lengths of the external branches is never greater than $2$, which means that it is largest for Kingman's coalescent.  The second inequality gives a rather sharp bound on the difference.  Recall that in Example \ref{exm1}, we have $\alpha_b \leq s \alpha$, so
$E[2 - J_n] \leq C' (\log n)/n$ for some other constant $C'$.  In Example \ref{exm2}, 
(\ref{alphab}) gives $\alpha_b \leq c(1 + \log b) \leq c(1 + \log n)$, which implies
$E[2 - J_n] \leq C''(\log n)^2/n$ for some constant $C''$.  Thus, in these examples, the lengths of the external branches are affected very little by multiple mergers when the sample size is large.  The reason is that, in large samples, a lot of coalescence occurs very quickly, so most ancestral lines have merged with at least one other ancestral line before the first multiple merger takes place.

\begin{proof}[Proof of Lemma \ref{singlem}] 
We start by proving the first inequality in (\ref{Lneq}) by induction.  As before, let $K_n$ be the amount of time that the integer $1$ is in a singleton block.  We need to show that $E[K_n] \leq 2/n$ for all $n \geq 2$.  First, note that $E[K_2] = \lambda_2^{-1} \leq 1$.  Now, suppose for some $n \geq 3$, we have $E[K_j] \leq 2/j$ for $j = 2, \dots, n-1$, and consider $E[K_n]$.  Let $T_n$ be the time of the first merger when the coalescent starts with $n$ blocks, and let $B \geq 2$ be the number of blocks involved in the merger at time $T_n$.  Note that $B$ is independent of $T_n$.  
Conditional on $B$, the probability that $1$ merges with at least one other block at time $T_n$ is $B/n$.  If this does not happen, then at least $n - B + 1$ blocks remain after the merger, so by the induction hypothesis, the expected time after $T_n$ that $\{1\}$ will remain a singleton is at most 
$2/(n-B+1)$.  Therefore,
$$E[K_n|T_n, B] \leq \bigg( \frac{B}{n} \bigg) T_n + \bigg( \frac{n-B}{n} \bigg) \bigg(T_n + 
\frac{2}{n-B+1} \bigg) = T_n + \frac{2(n-B)}{n(n-B+1)}.$$
Since $2 \leq B \leq n$, we have $(n-B)/(n-B+1) \leq (n-2)/(n-1)$.  Also,
$E[T_n] = \lambda_n^{-1} \leq 2/[n(n-1)]$, so
$$E[K_n]  \leq \frac{2}{n(n-1)} + \frac{2(n-2)}{n(n-1)} = \frac{2}{n},$$
which proves the first inequality.

The proof of the second inequality requires a coupling argument.  Let 
$(\Pi_n(t), t \geq 0)$ be the coalescent process defined in the statement of Proposition \ref{fuliprop}, and let $(\Upsilon_n(t), t \geq 0)$ be Kingman's coalescent, started from the partition of $1, \dots, n$ into singletons.  We may assume that the coalescent processes $\Pi_n$ and $\Upsilon_n$ are constructed
from Poisson processes $N_1$ and $N_2$ respectively on $(0, \infty) \times {\cal P}_n$, as described in section 3.  That is, whenever $(t, \pi)$ is a point of $N_1$, the partition $\Pi_n(t)$ is the coagulation of $\Pi_n(t-)$ by $\pi$, and whenever $(t, \pi)$ is a point of $N_2$, the partition $\Upsilon_n(t)$ is the coagulation of $\Upsilon_n(t-)$ by $\pi$.  Furthermore, these are the only jump times of $\Pi_n$ and $\Upsilon_n$.  Let $L_1$ and $L_2$ be the intensity measures of the second coordinate for the Poisson processes $N_1$ and $N_2$ respectively.  Then, for $\pi \in {\cal P}_n$, we have $L_2(\pi) = 1$ if $\pi$ consists of one block of size $2$ and $n-2$ singletons, and $L_2(\pi) = 0$ otherwise.  Also, $L_1(\pi) \geq L_2(\pi)$ for all $\pi \in {\cal P}_n$.  Therefore, we may assume that the Poisson processes $N_1$ and $N_2$ are coupled such that if $(t, \pi)$ is a point of $N_2$ then $(t, \pi)$ is a point of $N_1$.  The points $(t, \pi)$ in both $N_1$ and $N_2$ correspond to mergers in which two blocks coalesce at a time, while the points $(t, \pi)$ in $N_1$ but not $N_2$ correspond to multiple mergers caused by selective sweeps.

To compare the two processes,
note that $K_n = \inf\{t: \{1\} \mbox{ is not a singleton in }\Pi_n(t)\}$, and let
$K_n' = \inf\{t: \{1\} \mbox{ is not a singleton in }\Upsilon_n(t)\}$.  We have
$E[J_n] = n E[K_n]$.  By our previous results for Kingman's coalescent, we have
$E[K_n'] = 2/n$, and so $E[2 - J_n] = n E[K_n' - K_n]$.
Let $\tau = \inf\{t: \Pi_n(t) \neq \Upsilon_n(t)\}$, where we say $\tau = \infty$ if $\Pi_n(t) = \Upsilon_n(t)$ for all $t$.  For $\pi \in {\cal P}_n$, denote by $|\pi|$ the number of blocks in $\pi$.
Since $\Pi_n(t) = \Upsilon_n(t)$ for all $t \leq \tau$, we have
$$E[2 - J_n] = nE[K_n' - K_n] \leq nE[(K_n' - \tau)1_{\{\tau < K_n'\}}]
= n\sum_{b=2}^n E[(K_n' - \tau)1_{\{\tau < K_n'\}} 1_{\{|\Upsilon_n(\tau)| = b\}}].$$
For $b = 1, 2, \dots, n$, define $T_b = \inf\{t: |\Upsilon_n(t)| = b\}$.  If $\tau < K_n'$ and
$|\Upsilon_n(\tau)| = b$, then $K_n' > T_b$.  Therefore,
\begin{equation}
E[2 - J_n] \leq n \sum_{b=2}^n E[K_n' - \tau|\{\tau < K_n'\} \cap \{ |\Upsilon_n(\tau)| = b\}]
P(\{K_n' > T_b\} \cap \{ |\Upsilon_n(\tau)| = b\}).
\label{fueq1}
\end{equation}

If $\tau < K_n'$ and $|\Upsilon_n(\tau)| = b$, then $\{1\}$ is one of $b$ blocks of $\Upsilon_n(\tau)$, and by our previous results on Kingman's coalescent, the expected time before it merges with another block is $2/b$.  Thus, we have
\begin{equation}
E[K_n' - \tau|\{\tau < K_n'\} \cap \{ |\Upsilon_n(\tau)| = b\}] = \frac{2}{b}.
\label{fueq2}
\end{equation}
Note that $K_n' > T_b$ whenever $\{1\}$ remains a singleton at the time that Kingman's coalescent is down to $b$ blocks.  Whenever the coalescent goes from $j$ blocks to $j-1$, the probability that the integer $1$ is involved in the merger is $2/j$, so
\begin{equation}
P(K_n' > T_b) = \prod_{j=b+1}^n \bigg(1 - \frac{2}{j} \bigg) \leq
\exp \bigg( - \sum_{j=b+1}^n \frac{2}{j} \bigg) \leq
\exp \bigg( 1 - 2 \int_b^n \frac{1}{x} \: dx \bigg) = e \bigg( \frac{b}{n} \bigg)^2.
\label{fueq3}
\end{equation}
If $|\Upsilon_n(\tau)| = b$, then both $\Pi_n$ and $\Upsilon_n$ have the same $b$ blocks at time $T_b$, but at time $\tau$ the process $\Pi_n$ has a transition but $\Upsilon_n$ does not.  Since the total merger rate for $\Pi_n$ after time $T_b$ is $\lambda_b = \alpha_b + \binom{b}{2}$ and the total merger rate for $\Upsilon_n$ after time $T_b$ is $\binom{b}{2}$, we have
\begin{equation}
P(|\Upsilon_n(\tau)| = b|K_n' > T_b) \leq \frac{\alpha_b}{\lambda_b} \leq \frac{2 \alpha_b}{b(b-1)}.
\label{fueq4}
\end{equation}
Combining (\ref{fueq1})-(\ref{fueq4}), we get
$$E[2 - J_n] \leq n \sum_{b=2}^n \frac{4e \alpha_b b^2}{b^2(b-1)n^2} \leq \frac{C}{n} \sum_{b=2}^n
\frac{\alpha_b}{b},$$
which is the second inequality in (\ref{Lneq}).
\end{proof}

\begin{proof}[Proof of Proposition \ref{fuliprop}]  We have
\begin{equation}
E[S_n - h_{n-1} \eta_e] = (E[S_n] - \theta h_{n-1}) + h_{n-1}(\theta - E[\eta_e]) =
(E[S_n] - \theta h_{n-1}) + \frac{h_{n-1} \theta}{2} E[2 - J_n].
\label{Seta}
\end{equation}
By Proposition \ref{segprop}, $\lim_{n \rightarrow \infty} (E[S_n] - \theta h_{n-1}) = -\rho$.
It thus remains only to show that the second term on the right-hand side of (\ref{Seta}) goes to zero as $n \rightarrow \infty$.  Let $\epsilon > 0$.  By (\ref{segcond}), there exists a positive integer $N$ such that
$$\sum_{b=N}^{\infty} \frac{\alpha_b (1 + \log b)}{b^2} < \epsilon.$$  Therefore, by Lemma \ref{singlem},
\begin{align}
\limsup_{n \rightarrow \infty} \frac{h_{n-1} \theta}{2} E[2 - J_n] &\leq
\limsup_{n \rightarrow \infty} \frac{C h_{n-1} \theta}{2n} \sum_{b=2}^n \frac{\alpha_b}{b}
= \limsup_{n \rightarrow \infty} \frac{C h_{n-1} \theta}{2n} \bigg( \sum_{b=2}^N \frac{\alpha_b}{b}
+ \sum_{b=N}^n \frac{\alpha_b}{b} \bigg) \nonumber \\
&\leq 0 + \frac{C \theta}{2} \limsup_{n \rightarrow \infty} \sum_{b=N}^n \frac{\alpha_b h_{n-1}}{bn}
\leq \frac{C \theta}{2} \limsup_{n \rightarrow \infty} \sum_{b=N}^n \frac{\alpha_b (1 + \log b)}{b^2}
\leq \frac{C \theta \epsilon}{2}. \nonumber
\end{align}
Since this is true for all $\epsilon > 0$, and since $E[2 - J_n] \geq 0$ for all $n$ by
Lemma \ref{singlem}, we have $$\lim_{n \rightarrow \infty} \frac{h_{n-1} \theta}{2} E[2 - J_n] = 0,$$
which completes the proof of the proposition.
\end{proof}

We conclude this section with some comments about the power of Tajima's $D$-statistic and Fu and Li's $D$-statistic to detect selective sweeps.  The numerators of these two statistics, which are $\Delta_n - S_n/h_{n-1}$ and $S_n - h_{n-1} \eta_e$, each have mean zero when the ancestral process is Kingman's coalescent.  The expected values of these two numerators both converge to a negative constant as the sample size goes to infinity when multiple mergers can occur.  These statistics are used to test for departures from Kingman's coalescent.  If the goal is to test for multiple mergers caused by selective sweeps, one would reject the null hypothesis of no selective sweeps if the value of the statistic is too small (i.e. more negative than would be expected with Kingman's coalescent).

A natural question, then, is how much power these tests have to detect selective sweeps.  While a full analysis of this question would require a simulation study, we can obtain some insight from the analytical results presented above.  From the values of $a_n$ and $b_n$ in (\ref{tajeq}), which can be found in section 4.1 of Durrett (2002), we see that the standard deviation of the numerator of Tajima's $D$-statistic is $O(1)$ when the genealogy is given by Kingman's coalescent.  However, from the values of $c_n$ and $d_n$ in (\ref{fulieq}), which can be found in section 4.2 of Durrett (2002), we see that the numerator of Fu and Li's $D$-statistic has a standard deviation which is $O(\log n)$.  This means that, for large $n$, moderate negative values for the numerator of Fu and Li's $D$-statistic are not strong evidence against the null model of Kingman's coalescent, and thus a test based on Fu and Li's $D$-statistic will most likely have low power.
These observations are consistent with simulation results of Simonsen, Churchill, and Aquadro (1995), who found that Tajima's $D$-statistic has more power to detect selective sweeps than Fu and Li's $D$-statistic.

Neither of these tests has the desirable feature of many tests in classical statistics, which is that for all $\alpha > 0$, the power of the level $\alpha$ test tends to $1$ as the sample size $n$ tends to infinity.  Indeed, for the problem of detecting recurrent selective sweeps, no such test based on the genealogy of the sample can exist because, with positive probability, none of the selective sweeps affects the genealogy of the $n$ sampled lineages before we get back to the most recent common ancestor.  We formulate this observation precisely in the following proposition, which uses the coupling in the proof of Lemma \ref{singlem}.

\begin{Prop}
Let $(\Pi_n(t), t \geq 0)$ be the $\Lambda$-coalescent or $\Xi$-coalescent defined in the proof of Proposition \ref{fuliprop}, and assume that 
\begin{equation}
\sum_{b=2}^{\infty} \frac{\alpha_b}{b^2} < \infty,
\label{weakcond}
\end{equation}
which is slightly weaker than (\ref{segcond}).
Let $(\Upsilon_n(t), t \geq 0)$ be Kingman's coalescent, coupled with
$(\Pi_n(t), t \geq 0)$ as in the proof of Lemma \ref{singlem}.  Then there exists a constant $C > 0$ such that for all $n$, we have $P(\Upsilon_n(t) = \Pi_n(t) \mbox{ for all }t) \geq C$.  
\end{Prop}

\begin{proof}
Let $T_b = \inf\{t: |\Upsilon_n(t)| = b\}$.  Conditional on $\Pi_n(T_b) = \Upsilon_n(T_b)$, the probability that $\Pi_n(t) \neq \Upsilon_n(t)$ for some $t \in [T_b, T_{b-1}]$ is $\alpha_b/\lambda_b$.  It follows that
$$P(\Upsilon_n(t) = \Pi_n(t) \mbox{ for all }t) = \prod_{b=2}^n \bigg(1 - \frac{\alpha_b}{\lambda_b} \bigg).$$  Note that $\alpha_b/\lambda_b \leq 2 \alpha_b/[b(b-1)]$ for all $b$.
By (\ref{weakcond}), there exists a positive integer $N$ such that $6 \alpha_b/[b(b-1)] \leq 1$ for all $b \geq N$, and if $0 \leq x \leq 1$ then $1 - x/3 \geq e^{-x}$.  Putting these results together, we get
\begin{align}
P(\Upsilon_n(t) = \Pi_n(t) \mbox{ for all }t) &\geq \prod_{b=2}^{N-1} \bigg(1 - \frac{\alpha_b}{\lambda_b} \bigg) \prod_{b=N}^{\infty} \exp \bigg( - \frac{6 \alpha_b}{b(b-1)} \bigg) \nonumber \\
&= \prod_{b=2}^{N-1} \bigg(1 - \frac{\alpha_b}{\lambda_b} \bigg)
\exp \bigg(- \sum_{b=N}^{\infty} \frac{6 \alpha_b}{b(b-1)} \bigg) \geq C,
\nonumber
\end{align}
where the last inequality uses (\ref{weakcond}) again.
\end{proof}

\section{Proofs of convergence theorems}

In this section, we prove Theorem \ref{mainth} and Proposition \ref{coalprop}.
The proofs use Propositions \ref{sweepprop1} and \ref{sweepprop2} in combination with the Poisson
process construction of coalescents with multiple or simultaneous multiple
collisions.  

Recall the model presented in subsection 2.1 of how the population behaves
following a single beneficial mutation.  As in subsection 2.1, assume for
now that a beneficial mutation occurs at time $0$.  Let $X(t)$ be the
number of chromosomes with the favorable $B$ allele at time $t$, and let
$\tau = \inf\{t: X(t) \in \{0, 2N\}\}$.  Let
$0 = \xi_0 < \xi_1 < \xi_2 < \dots$
be the times of the proposed replacements, which occur at times of a rate
$2N$ Poisson process.  Let $0 = \xi_0'
< \xi_1' < \xi_2' < \dots$ be the subset of
these times at which the number of individuals with the favorable allele
changes.  As observed in Schweinsberg and Durrett (2004), if $1 \leq k
\leq 2N-1$, then $P(X(\xi_{i+1}') = k+1|X(\xi_i') = k) = 1/(2-s)$ and
$P(X(\xi_{i+1}') = k-1|X(\xi_i') = k) = (1-s)/(2-s)$.  Thus, the number
of chromosomes with the $B$ allele behaves like an asymmetric random walk
until it reaches $0$ or $2N$.  For integers $i$, $j$, and $k$
such that $0 \leq i \leq k \leq j \leq 2N$ and $i < j$, define
$$p(i,j,k) = P(\inf\{s \geq t: X(s) = j\} < \inf\{s \geq t: X(s) = i\}|X(t) = k),$$
which is the probability that if at some time there are $k$ chromosomes
with $B$, the number of $B$'s will reach $j$ before $i$.  Using
the fact that $(1-s)^{\xi_n'}$ is a martingale and applying
the Optional Sampling Theorem, we get (see also Durrett (2002) or
Lemma 3.1 of Schweinsberg and Durrett (2004))
$$p(i,j,k) = \frac{1 - (1-s)^{k-i}}{1 - (1-s)^{j-i}}.$$  Therefore, the
probability that the beneficial mutation leads to a selective sweep is
$p(0,2N,1) = s/(1 - (1-s)^{2N})$.

Lemma \ref{taulem} below shows that the length of time that the beneficial
allele is present in the population is only $O(\log N)$.  Since we speed
up time by a factor of $N$ to define the ancestral process, it will follow
that for large populations, on the time scale of the ancestral process
the lineages that coalesce as a result of a selective sweep coalesce almost
at the same time.  It is well-known
(see Durrett (2002)) that a selective sweep takes time approximately
$(2/s) \log(2N)$.  However, since a beneficial mutation leads to a
selective sweep with probability approximately $s$, we get a bound on
$E[\tau]$ that does not depend on $s$.

\begin{Lemma}
We have $E[\tau] \leq 4( \log N + 1)$.
\label{taulem}
\end{Lemma}

\begin{proof}
For $1 \leq k \leq 2N - 1$, let $S_k = \# \{i \geq 0: X(\xi_i') = k\}$ and
$T_k = \# \{i \geq 0: X(\xi_i) = k\}$, where $\# S$ denotes the cardinality of
a set $S$.  Let $q_k = P(X(\xi_j') \neq k \mbox{ for all }j > i|X(\xi_i') = k)$
be the probability that the asymmetric random walk never returns to $k$.
Note that $E[S_k|S_k \geq 1] = 1/q_k$.
Also, $P(X(\xi_i) = k \mbox{ for some }k) = p(0,k,1) = s/(1 - (1-s)^k)$.
Therefore,
\begin{equation}
E[S_k] = P(S_k \geq 1) E[S_k|S_k \geq 1] =
\frac{s}{q_k(1 - (1-s)^k)}.
\label{ESk}
\end{equation}
We have, for $1 \leq k \leq 2N-1$,
\begin{align}
q_k &= \bigg( \frac{1-s}{2-s} \bigg) [1 - p(0,k,k-1)] +
\bigg( \frac{1}{2-s} \bigg) p(k, 2N, k+1) \nonumber \\
&= \bigg( \frac{1-s}{2-s} \bigg) \bigg[ 1 - \frac{1 - (1-s)^{k-1}}{1 - (1-s)^k}
\bigg] + \bigg( \frac{1}{2-s} \bigg) \frac{1 - (1-s)}{1 - (1-s)^{2N-k}}
\nonumber \\
&= \bigg( \frac{1-s}{2-s} \bigg) \frac{s(1-s)^{k-1}}{1 - (1-s)^k} + 
\bigg( \frac{1}{2-s} \bigg) \frac{s}{1 - (1-s)^{2N-k}} \nonumber \\
&\geq \frac{s}{2-s} \bigg( \frac{(1-s)^k}{1 - (1-s)^k} + 1 \bigg)
= \frac{s}{(2-s)(1 - (1-s)^k)}. \nonumber
\end{align}
It follows from this result and (\ref{ESk}) that $E[S_k] \leq 2-s$
for all $k$.
Schweinsberg and Durrett (2004) calculated that $P(X(\xi_{i+1}) \neq
X(\xi_i)|X(\xi_i) = k) = k(2N-k)(2-s)/(2N)^2$.  It follows that
$$E[T_k] = E[S_k] \bigg( \frac{(2N)^2}{k(2N-k)(2-s)} \bigg) \leq
\frac{4N^2}{k(2N-k)}.$$  Since $E[\xi_{i+1} - \xi_i] = 1/2N$ for all $i$, we
have $$E[\tau] = \frac{1}{2N} \sum_{k=1}^{2N-1} E[T_k] \leq \sum_{k=1}^{2N-1}
\frac{2N}{k(2N-k)} \leq 2 \sum_{k=1}^N \frac{2}{k} \leq 4(\log N + 1),$$
as claimed.
\end{proof}

We now use this result to prove part 2 of Proposition \ref{sweepprop1},
which shows that beneficial mutations do not cause lineages to coalesce
when the beneficial gene dies out.

\begin{proof}[Proof of part 2 of Proposition \ref{sweepprop1}]
Suppose $X(\tau) = 0$ and $\Theta \neq \kappa_0$.  Then it can not be true
that for all $t \in [0, \tau]$, the $n$ individuals sampled at time $\tau$
all have distinct ancestors with the $b$-chromosome at time $t$.  Therefore,
there is an
integer $i$ with $\xi_i \leq \tau$ such that one of the following is true:
\begin{enumerate}
\item The ancestor at time $\xi_i$ of one of the individuals sampled at
time $\tau$ has the $b$ allele, but the ancestor of the same individual
at time $\xi_{i-1}$ has the $B$ allele because of recombination.

\item There are two individuals in the sample at time $\tau$ that have
distinct ancestors with the $b$ allele at time $\xi_i$, but both of them
have the same ancestor at time $\xi_{i-1}$.
\end{enumerate}

We now calculate the probability of these events conditional on
$X(\xi_i) = k$, where $1 \leq k \leq N^{1/2}$.  We assume $N \geq 2$.
For a randomly chosen $b$ chromosome at time $\xi_i$ to have a $B$
chromosome as its ancestor at time $\xi_{i-1}$, the chosen $b$ chromosome
must be the new one born at time $\xi_i$ (which has probability at most
$1/(2N-k)$ because $2N-k$ chromosomes have the $b$ allele at time
$\xi_i$), there must be recombination at this
time (which happens with probability $r$), and the ancestor at the site
of interest must be a $B$ chromosome (which happens with probability at
most $(k+1)/2N$ because $X(\xi_{i-1}) \leq k+1$).  Therefore, the probability
that all three events occur is at most
$r(k+1)/[(2N-k)(2N)] \leq r/N^{3/2}$.  Also, at most one pair of $b$
chromosomes at time $\xi_i$ can have the same ancestor at time $\xi_{i-1}$,
so the probability that two randomly chosen $b$ chromosomes coalesce at this
time is at most $\binom{2N-k}{2}^{-1} = 2/[(2N-k)(2N-k-1)] \leq 2/N^2$.

By Lemma \ref{taulem}, if $M$ is the integer such that $\xi_M = \tau$,
then $E[M] \leq (2N)[4(\log N + 1)] = 8N(\log N + 1)$.  Since there
are $n$ individuals and $\binom{n}{2}$ pairs in the sample, combining
these bounds gives
\begin{equation}
P(X(\tau) = 0, X(t) \leq N^{1/2} \mbox{ for all }t,
\mbox{ and } \Theta \neq \kappa_0)
\leq 8N(\log N + 1) \bigg( \frac{nr}{N^{3/2}} + \frac{n(n-1)}{N^2} \bigg).
\label{nosweep1}
\end{equation}
Note that for $1 \leq k \leq 2N-1$, we have
\begin{align}
P(X(\tau) = 0 \mbox{ and }X(t) = k \mbox{ for some }t) &\leq
P(X(\tau) = 0|X(t) = k \mbox{ for some }t) \nonumber \\
&= 1 - p(0, 2N,k) = 1 - \frac{1 - (1-s)^k}{1 - (1-s)^{2N}} \leq (1-s)^k.
\nonumber
\end{align}
Therefore,
\begin{equation}
P(X(\tau) = 0 \mbox{ and }
X(t) > N^{1/2} \mbox{ for some }t) \leq (1-s)^{N^{1/2}}.
\label{nosweep2}
\end{equation}
Combining (\ref{nosweep1}) and (\ref{nosweep2}), we get
$$P(X_{\tau} = 0 \mbox{ and } \Theta \neq \kappa_0) \leq
(1 - s)^{N^{1/2}} + 8N(\log N + 1) \bigg( \frac{nr}{N^{3/2}} +
\frac{n(n-1)}{N^2} \bigg).$$
Part 2 of Proposition \ref{sweepprop1} follows because 
$r \leq C' \log(2N)$ and $s$ is fixed.
\end{proof}

We now consider our model of recurrent selective sweeps and work towards
the proof of Theorem \ref{mainth}.  We will first define a coalescent
with multiple collisions.  We will then show that this process can be
coupled with the ancestral process $(\Psi_N(t), t \geq 0)$ such that,
given a finite number of times $0 < u_1 < \dots < u_m$, the processes
agree at these times with high probability.

Recall that $K_N$ is a Poisson point process on $\R \times [-L, L]
\times [0,1]$ with intensity $\lambda \times \mu_N$.  We can define another
Poisson point process $K_N^*$ on $[0, \infty) \times [-L, L] \times [0,1]$
which consists of all the points $(-t/N, x, s)$ such that $(t, x, s)$ is
a point of $K_N$ and $t \leq 0$.  By the Mapping Theorem for Poisson processes
(see section 2.3 of Kingman (1993)), $K_N^*$ is a Poisson process with
intensity measure $\lambda \times N \mu_N$.  The points in
$K_N^*$ can be ordered by their first coordinate, so we can write the
points as $(t_i, x_i, s_i)$ for positive integers $i$, where
$0 < t_1 < t_2 < \dots$ a.s.  Also, define $t_0 = 0$.

We now define a ${\cal P}_n$-valued coalescent process
$\Pi_N = (\Pi_N(t), t \geq 0)$.  Let $\Pi_N(0)$ be the partition $\kappa_0$ 
of $\{1, \dots, n\}$ into singletons.  Given $\Pi_N(t_i)$ for some $i \geq 0$,
we define $\Pi_N(t)$ for $t_i < t \leq t_{i+1}$ in two steps.  First,
we let the process obey the law of Kingman's coalescent
over the interval $(t_i, t_{i+1})$, meaning that each possible transition
that involves the merging of two blocks happens at rate one.
Second, let $\pi_{i+1}$ be a random partition of $\{1, \dots, n\}$,
independent of $(\Pi_N(t), 0 \leq t < t_{i+1})$, such that for an event
$A_{i+1}$ of probability $s_{i+1}$, we have $\pi_{i+1} = \kappa_0$
on $A_{i+1}^c$ and the conditional distribution of $\pi_{i+1}$ given
$A_{i+1}$ is $Q_{p,n}$, where $p = e^{-r_N(x_{i+1}) \log(2N)/s_{i+1}}$.
We then define $\Pi_N(t_{i+1})$ to be the coagulation of
$\Pi_N(t_{i+1}-)$ by $\pi_{i+1}$.

The lemma below states that the coalescent process $\Pi_N$ that we have
just defined is a coalescent with multiple collisions.

\begin{Lemma}
Let $\eta_N$ be the measure on $(0,1]$ such that
$$\eta_N([y,1]) = \int_{-L}^L \int_0^1 
s 1_{\{e^{-r_N(x) \log(2N)/s} \geq y\}} \: N \mu_N(dx \times ds)$$
for all $y \in (0,1]$.  Let $\Lambda_{0,N}$ be the measure on
$(0, 1]$ such that $\Lambda_{0,N}(dx) = x^2 \eta_N(dx)$, and
let $\Lambda_N = \delta_0 + \Lambda_{0,N}$.  Then the process
$(\Pi_N(t), t \geq 0)$ is the  ${\cal P}_n$-valued $\Lambda_N$-coalescent.
\label{coallem}
\end{Lemma}

\begin{proof}
Let $K_N'$ be the point process on $[0, \infty) \times {\cal P}_n$
consisting of the points $(t_i, \pi_i)$.
By the Marking Theorem for Poisson processes (see section 5.2 of
Kingman (1993)), $K_N'$ is also a Poisson point process.  Given
$(t_i, x_i, s_i)$, the partition $\pi_i$ has distribution
$Q_{p,n}$, where $p = e^{-r_N(x_i) \log(2N)/s_i}$, conditional on an
event of probability $s_i$ and otherwise is $\kappa_0$.
Therefore, the intensity measure of $K_N'$ is given by $\lambda \times H$,
where, for $\pi \neq \kappa_0$,
\begin{align}
H(\pi) &= \int_{-L}^L \int_0^1 
s Q_{e^{-r_N(x) \log(2N)/s}, n}(\pi) \: N \mu_N(dx \times ds) \nonumber \\
&= \int_0^1 Q_{p,n}(\pi) \: \eta_N(dp) = \int_0^1 Q_{p,n}(\pi)
p^{-2} \Lambda_{0,N}(dp). \nonumber
\end{align}
By comparing this with (\ref{cmcpois}) and recalling that
$\Pi_N$ follows the law of Kingman's coalescent over the intervals
$(t_{i-1}, t_i)$, we conclude that $\Pi_N$ is the $\Lambda_N$-coalescent.
\end{proof}

The next lemma states that it is unlikely for there to be a beneficial
allele in the population at any fixed time.
Recall that ${\cal T}_N = \{t: (t,x,s) \mbox{ is a point in }K_N
\mbox{ for some } x \mbox{ and }s\}$.

\begin{Lemma}
There exists a constant $C$, not depending on $N$, such that for any
fixed $y \in \R$, we have $P(y \in [t, \tau(t)] \mbox { for some }t \in
{\cal T}_N) \leq (C \log N)/N.$
\label{overlap}
\end{Lemma}

\begin{proof}
The points of ${\cal T}_N$ form a Poisson process on $\R$ of rate
$\gamma_N$, where $\gamma_N = \mu_N([-L,L] \times [0,1])$.
Recall from Lemma \ref{taulem} that if $\tau$ denotes the amount of time for
which a beneficial allele is present in between $1$ and $2N-1$ members of
the population, then $E[\tau] \leq 4(\log N + 1)$.  Therefore,
\begin{align}
P(y \in [t, \tau(t)] \mbox{ for some }t \in {\cal T}_N) &\leq
\int_{-\infty}^y P(\tau \geq y - x) \gamma_N \: dx \nonumber \\
&= \gamma_N \int_0^{\infty} P(\tau \geq x) \: dx = \gamma_N E[\tau]
\leq 4 \gamma_N (\log N + 1). \nonumber
\end{align}
Since the measures $N \mu_N$ converge weakly to $\mu$, the sequence
$(N \gamma_N)_{N=1}^{\infty}$ converges to $\mu([-L,L] \times [0,1])$
and therefore is bounded.  The lemma follows.
\end{proof}

We now show how to couple the processes $\Psi_N$ and $\Pi_N$ so that
they agree at a given finite set of times with high probability.  We first
consider how the ancestral process $\Psi_N$ behaves around the times
$t_1, t_2, \dots$.  For positive integers $i$, let $\tau_i = -\tau(-Nt_i)/N$.
We have $-Nt_i \in {\cal T}_N$.  However, recall from subsection 2.2 that
not all points in ${\cal T}_N$ are in ${\cal T}'_N$ because some potential
mutations are discarded to avoid overlapping selective sweeps. 
When  $-Nt_i \in {\cal T}_N'$, there is a beneficial allele in the
population during the time interval $[-Nt_i, \tau(-Nt_i))$, and this
affects the process $\Psi_N$ over the interval $[\tau_i, t_i]$.

For each $i$ such that $-Nt_i \in {\cal T}_N'$, we can define a random
partition $\theta_i \in {\cal P}_n$ by choosing $n$ individuals from the
population at time $\tau(-Nt_i)$ and declaring two integers $j$ and $k$
to be in the same block of $\theta_i$ if and only if the $j$th and
$k$th individuals chosen got their
allele at the neutral site of interest from the same ancestor at
time $-Nt_i$.  If $\tau_i > 0$ and the partition $\Psi_N(\tau_i)$ contains
$b_i$ blocks, we can choose the $n$ individuals at time $\tau(-Nt_i)$
by first picking the $b_i$ individuals that are ancestors of the $n$
individuals that were sampled at time zero, and then choosing the
remaining $n - b_i$ at random.  This will ensure that, for $i$ such that 
$-Nt_i \in {\cal T}_N'$ and $\tau_i > 0$, the random partition
$\Psi_N(t_i)$ is the coagulation of $\Psi_N(\tau_i)$ by $\theta_i$.

Moreover, the conditional distribution of $\theta_i$
given $(t_i, x_i, s_i)$ and given that $-Nt_i \in {\cal T}_N'$
is the same as the distribution of the random partition $\Theta$
defined in subsection 2.1, when the selective advantage is $s_i$ and
the recombination probability is $r_N(x_i)$.  Recall that when a beneficial
mutation occurs in the population with selective advantage $s_i$, it spreads
to the entire population with probability $s_i/(1 - (1 - s_i)^{2N})$.
Therefore, by Proposition \ref{sweepprop1}, 
the distribution of $\Theta$ is approximately that of a random partition
that has distribution $Q_{p,n}$, where $p = e^{-r_N(x_i) \log(2N)/s_i}$,
on an event of probability $s_i/(1 - (1 - s_i)^{2N})$ and is $\kappa_0$
on the complementary event.  However, this is the same as the conditional
distribution of $\pi_i$ given $(t_i, x_i, s_i)$, except we have
$s_i/(1 - (1 - s_i)^{2N})$ instead of $s_i$.  It thus follows
from Proposition \ref{sweepprop1} that
we can couple the partitions $\theta_i$ and $\pi_i$ such that
for any $\delta > 0$,
\begin{equation}
P(\theta_i \neq \pi_i \mbox{ and }-Nt_i \in {\cal T}_N'|(t_i, x_i, s_i))
\leq \frac{C_{\delta}}{\log N} + 1_{\{s_i < \delta\}},
\label{couppart}
\end{equation}
where $C_{\delta}$ is a constant that depends on $\delta$.  Note that
we only get the $O(1/(\log N))$ bound when $s_i \geq \delta$ because
of the assumption in Proposition \ref{sweepprop1} that $s$ is fixed.

Finally, we consider the processes during the intervals
$(t_i, t_{i+1})$.  The process $\Pi_N$ behaves like Kingman's coalescent
during these intervals.
Let $${\cal I}_N^* = \bigcup_{i=1}^{\infty} [\tau_i, t_i].$$
The process $\Psi_N$ behaves like Kingman's coalescent during the
intervals in $(0, \infty) \setminus {\cal I}_N^*$ because the population
follows the Moran model during the corresponding intervals.
Therefore, if $\Pi_N(t_i) = \Psi_N(t_i)$, we can couple the processes
so that $\Pi_N(t) = \Psi_N(t)$ for all $t \in [t_i, \phi_i)$,
where $\phi_i = \inf\{t > t_i: t \in {\cal I}_N^* \}$.

\begin{Prop}
Suppose the processes $\Pi_N$ and $\Psi_N$ are coupled in the manner
described above.  Let $0 < u_1 < \dots < u_m$ be fixed times.
Let $\epsilon > 0$.  For sufficiently large $N$, we have
\begin{equation}
P(\Pi_N(u_i) \neq \Psi_N(u_i) \mbox{ for some }i \in \{1, \dots, m\}) 
< \epsilon.
\label{maincoup}
\end{equation}
\label{coupprop}
\end{Prop}

\begin{proof}
Let $K = \sup\{k: t_k \leq u_m\}$.  Suppose the following conditions hold:
\begin{enumerate}
\item For $i = 1, \dots, m$, we have $u_i \notin {\cal I}_N^*$.

\item For all positive integers $i$, we have $\tau_i > 0$.

\item For $i = 1, \dots, K$, we have $-Nt_i \in {\cal T}_N'$.

\item For $i = 1, \dots, K$, we have $\Pi_N(\tau_i) = \Pi_N(t_i-)$.

\item For $i = 1, \dots, K$, we have $\theta_i = \pi_i$.
\end{enumerate}
Conditions 2 and 3 imply that
$$0 = t_0 < \tau_1 < t_1 < \tau_2 < t_2 < \dots < \tau_K < t_K \leq u_m.$$
Condition 1 with $i = m$ implies further that $\tau_j > u_m$ for all $j > K$,
so $(t_K, u_m] \subset \R \setminus {\cal I}_N^*$.
We know that $\Pi_N(t_0) = \Psi_N(t_0) = \kappa_0$.  Suppose, for
some $i \in \{0, \dots, K-1\}$, that $\Pi_N(t_i) = \Psi_N(t_i)$.
Then the coupling gives $\Pi_N(t) = \Psi_N(t)$ for all
$t \in [t_i, \tau_{i+1})$.  Condition 4 gives
$\Pi_N(\tau_{i+1}) = \Pi_N(t_{i+1}-)$.  Conditions 2 and 3 imply that
$\Psi_N(t_{i+1})$ is the coagulation of $\Psi_N(\tau_{i+1})$ by $\theta_{i+1}$.
Since $\Pi_N(t_{i+1})$ is the coagulation of $\Pi_N(t_{i+1}-)$ by $\pi_{i+1}$,
condition 5 ensures that $\Pi_N(t_{i+1}) = \Psi_N(t_{i+1})$.
Thus, $\Pi_N(t_i) = \Psi_N(t_i)$ for $i = 0, 1, \dots, K$, 
and the coupling combined with the fact that $\Pi_N(t_K) = \Psi_N(t_K)$
gives $\Pi_N(t) = \Psi_N(t)$ for all $t \in (t_K, u_m]$.  Thus,
we have $\Pi_N(t) = \Psi_N(t)$ for all $t \in [0, u_m] \setminus {\cal I}_N^*$.
Therefore, by condition 1, $\Pi_N(u_i) = \Psi_N(u_i)$ for $i = 1, \dots, m$.
It thus remains only to show that conditions 1 through 5 occur with
high probability.  For the rest of the proof, we allow the constant $C$
to change from line to line.

If $u_i \in {\cal I}_N^*$, then there exists $t \in {\cal T}_N$ such that
$-Nu_i \in [t, \tau(t)]$.  Therefore, by Lemma \ref{overlap},
$$P(u_i \in {\cal I}_N^* \mbox{ for some }i \in \{1, \dots, m\}) \leq
\frac{C \log N}{N}.$$  Likewise,
if $\tau_i < 0$ for some $i$, then $-Nt_i \leq 0 \leq \tau(-Nt_i)$
and $-Nt_i \in {\cal T}_N$.  It follows that $P(\tau_i < 0 \mbox{ for some }i) 
\leq C(\log N)/N$ by Lemma \ref{overlap}.

To deal with conditions 3, 4, and 5, let $l_j = ju_m/N$ for
$j = 0, 1, \dots, N$, and define the intervals $I_1, \dots, I_N$
by $I_j = [l_{j-1}, l_j]$.  Note that the number of the points $t_i$
in an interval $I_j$ is Poisson with mean $u_m \gamma_N$.  Therefore,
the probability that some point $t_i$ falls in $I_j$ is at most
$u_m \gamma_N \leq C/N$.
The probability that two or more points fall in $I_j$ is at most
$u_m^2 \gamma_N^2 \leq C/N^2$.
If there is a point $t_i \in I_j$ with $-Nt_i \notin {\cal T}_N'$,
then either there are two points in $I_j$ or there is one point in
$I_j$ and $l_j \in {\cal I}_N^*$.  The event that there is at least
one point in $I_j$ is independent of the event that $l_j \in {\cal I}_N^*$,
so using Lemma \ref{overlap} again, the probability that both occur
is at most $C(\log N)/N^2$.

When the number of blocks in the coalescent is at most $n$, the
total transition rate of the process $\Pi_N$ is bounded by
$\binom{n}{2} + N \gamma_N$.  The probability that there is any
point $t_i$ in $I_j$ is at most $C/N$, so by Lemma \ref{taulem},
\begin{equation}
P(\Pi_N(\tau_i) \neq \Pi_N(t_i-) \mbox{ for some }t_i \in I_j) \leq 
\frac{C}{N} \bigg(
\binom{n}{2} + N\gamma_N \bigg) E[t_i - \tau_i] \leq \frac{C \log N}{N^2}.
\nonumber
\end{equation}
Finally, we may choose $\delta$ small enough that
$P(s_i \leq \delta) < \epsilon$, and then (\ref{couppart}) gives
$$P(\theta_i \neq \pi_i|t_i \in I_j) < \epsilon + \frac{C_{\delta}}{\log N},$$
where $C_{\delta}$ is a constant that depends on $\delta$.  Therefore,
$P(\theta_i \neq \pi_i \mbox{ for some } t_i \in I_j)
\leq \epsilon/N + C_{\delta}/(N \log N)$.  Since there are only
$N$ intervals $I_j$, we can add these bounds to show that the probability
that conditions 1 through 5 all hold is at least $1 - \epsilon$ for
sufficiently large $N$, which implies the statement of the proposition.
\end{proof}

\begin{proof}[Proof of Theorem \ref{mainth}]
Let $0 < u_1 < \dots < u_m$ be fixed times.  Let $\epsilon > 0$.
Define $\Lambda_N$ as in Lemma \ref{coallem}, and let $\Pi_N$ be
a ${\cal P}_n$-valued $\Lambda_N$-coalescent.  In view of Proposition
\ref{coupprop}, it suffices to show that for all $\pi_1, \dots, \pi_m 
\in {\cal P}_n$, we have
$$|P(\Pi_N(u_i) = \pi_i \mbox{ for all } i \in \{1, \dots, m\}) -
P(\Pi(u_i) = \pi_i \mbox{ for all } i \in \{1, \dots, m\})| < \epsilon$$
for sufficiently large $N$.  Therefore (see Pitman (1999)), it
suffices to show that the measures $\Lambda_N$ converge weakly to
$\Lambda$.  Thus, we need to show (see Billingsley (1999), Theorem 2.1)
that for any bounded uniformly continuous function $h$ on $[0,1]$,
we have $\int_0^1 h(x) \: \Lambda_N(dx) \rightarrow \int_0^1 h(x) \: \Lambda(dx)$
as $N \rightarrow \infty$.  By the definitions of $\Lambda_N$ and
$\Lambda$, it suffices to show that
$\int_0^1 h(x) \: \eta_N(dx) \rightarrow \int_0^1 h(x) \: \eta(dx)$
as $N \rightarrow \infty$  for
any bounded uniformly continuous function $h$ on $(0,1]$.
By the definitions of the measures
$\eta_N$ and $\eta$, this is equivalent to showing that
\begin{equation}
\lim_{N \rightarrow \infty} \int_{-L}^L \int_0^1 s h \big(
e^{-r_N(x) \log(2N)/s} \big) \: N \mu_N(dx \times ds) =
\int_{-L}^L \int_0^1 s h \big(e^{-r(x)/s} \big) \: \mu(dx \times ds)
\label{weakconv}
\end{equation}
for any bounded uniformly continuous function $h$ on $(0,1]$.
However, it is easy to deduce (\ref{weakconv}) from the
boundedness and uniform continuity of $h$, the uniform convergence of
$(\log 2N) r_N$ to $r$, the continuity of $r$, and the weak convergence
of the measures $N \mu_N$ to $\mu$.
\end{proof}

\begin{proof}[Proof of Proposition \ref{coalprop}]  Proposition
\ref{coalprop} can be proved by repeating the proof of Proposition
\ref{coupprop} with minor changes.  To prove the
first part of the proposition, we construct the coalescent process $\Pi_N$
as before.  Because $N \mu_N = \mu$ and $\log(2N) r_N = r$
for all $N$, we have $\Lambda_N = \Lambda$ for all $N$. 
It follows from Lemma \ref{coallem} that $\Pi_N$ is a $\Lambda$-coalescent
for all $N$.  Thus, it suffices to show (\ref{maincoup}), but with
$C/(\log N)$ on the right-hand side instead of $\epsilon$.
Because we are assuming that $\mu$ is concentrated on
$[-L, L] \times [\epsilon, 1]$ for some $\epsilon > 0$, we can
choose $\delta = \epsilon$ and drop the indicator from the
right-hand side of (\ref{couppart}) to get a bound of
$C_{\epsilon}/(\log N)$.  We then obtain $C/(\log N)$ on the
right-hand side of (\ref{maincoup}) by following the same steps as before.

To prove the second part of Proposition \ref{coalprop}, we modify the
definition of $\Pi_N$.  Conditional on $A_i$,
we give $\pi_i$ the distribution
$Q_{R(r_N(x_i)/s_i, \lfloor 2N s_i \rfloor), n}$.  We set $\pi_i = \kappa_0$
on $A_i^c$.
The intensity measure of $K_N'$ is then given by $\lambda \times J$, where,
for all $\pi \neq \kappa_0$, we have
\begin{align}
J(\pi) &= \int_{-L}^L \int_0^1 s
Q_{R(r_N(x)/s, \lfloor 2N s \rfloor), n}(\pi) \: N \mu_N(dx \times ds)
\nonumber \\
&= \int_{\Delta} Q_{\delta_x, n}(\pi) \: G_N(dx)
= \int_{\Delta} Q_{\delta_x, n}(\pi)
\bigg( \sum_{j=1}^{\infty} x_j^2 \bigg)^{-1} \: \Xi_{N,0}(dx). \nonumber
\end{align}
By comparing this with (\ref{csmcpois}),
we see that the process $\Pi_N$ is a $\Xi_N$-coalescent.
It follows from Proposition \ref{sweepprop2} that we
can replace $C_{\delta}/(\log N)$ on the right-hand side of (\ref{couppart})
by $C_{\delta}/(\log N)^2$.  This gives the second part of the proposition.
\end{proof}

\bigskip
\begin{center}
{\bf {\Large References}}
\end{center}

\mn N. H. Barton (1995).  Linkage and the limits to natural selection.
{\it Genetics}, {\bf 140}, 821-841.

\mn N. H. Barton (1998). The effect of hitch-hiking on neutral genealogies.
{\it Genet. Res.} {\bf 72}, 123-133.

\mn N. H. Barton, A. M. Etheridge, and A. K. Sturm (2004).  Coalescence in a
random background.  {\it Ann. Appl. Probab}. {\bf 14}, 754-785.

\mn J. Bertoin and J.-F. LeGall (2003).  Stochastic flows associated
to coalescent processes.  {\it Probab. Theory Relat. Fields},
{\bf 126}, 261-288.

\mn P. Billingsley (1999).  {\it Convergence of Probability Measures}.
2nd ed.  New York, Wiley.

\mn J. M. Braverman, R. R. Hudson, N. L. Kaplan, C. H. Langley,
and W. Stephan (1995).  The hitchhiking effect on the site frequency spectrum
of DNA polymorphisms.  {\it Genetics}, {\bf 140}, 783-796.

\mn P. Donnelly and T. G. Kurtz (1999).  Genealogical processes for
Fleming-Viot models with selection and recombination.  {\it Ann. Appl.
Prob.} {\bf 9}, 1091-1148.

\mn R. Durrett (2002).  {\it Probability Models for DNA Sequence
Evolution}.  New York, Springer-Verlag.

\mn R. Durrett and J. Schweinsberg (2004a).  Approximating selective
sweeps.  {\it Theor. Popul. Biol.} {\bf 66}, 129-138.

\mn R. Durrett and J. Schweinsberg (2004b).  Power laws for family sizes
in a duplication model.  Preprint. 
Available at http://front.math.ucdavis.edu/math.PR/0406216.

\mn Y. X. Fu and W. H. Li (1993).  Statistical tests of neutrality of mutations.
{\it Genetics}, {\bf 133}, 693-709.

\mn P. J. Gerrish and R. E. Lenski (1998).  The fate of competing
beneficial mutations in an asexual population.  {\it Genetica},
{\bf 102/103}, 127-144.

\mn J. H. Gillespie (2000).  Genetic drift in an infinite population:
the pseudohitchhiking model.  {\it Genetics}, {\bf 155}, 909-919.

\mn N. L. Kaplan, R. R. Hudson, and C. H. Langley (1989). 
The ``hitchhiking effect'' revisited.  {\it Genetics}, {\bf 123}, 887-899. 

\mn Y. Kim and W. Stephan (2003).  Selective sweeps in the presence
of interference among partially linked loci.  {\it Genetics},
{\bf 164}, 389-398.

\mn J. F. C. Kingman (1978).  The representation of partition
structures.  {\it J. London Math. Soc.} {\bf 18}, 374-380.

\mn J. F. C. Kingman (1982).  The coalescent.
{\it Stochastic Process. Appl.} {\bf 13}, 235-248.

\mn J. F. C. Kingman (1993).  {\it Poisson processes}. Oxford, Clarendon Press.

\mn S. M. Krone and C. Neuhauser (1997).  Ancestral processes with
selection.  {\it Theor. Popul. Biol.}  {\bf 51}, 210-237.

\mn J. Maynard Smith and J. Haigh (1974).  The hitchhiking effect of a
favorable gene.  {\it Genet. Res.} {\bf 23}, 23-35.

\mn M. M{\"o}hle and S. Sagitov (2001).
A classification of coalescent processes for haploid exchangeable
population models.  {\it Ann. Probab.} {\bf 29}, 1547-1562.

\mn P. A. P. Moran (1958).  Random processes in genetics.
{\it Proc. Cambridge Philos. Soc.} {\bf 54}, 60-71.

\mn C. Neuhauser and S. M. Krone (1997).  The genealogy of samples in
models with selection.  {\it Genetics}, {\bf 145}, 519-534.

\mn J. Pitman (1999).  Coalescents with multiple collisions.
{\it Ann. Probab.} {\bf 27}, 1870-1902.  

\mn M. Przeworski (2002).  The signature of positive selection at
randomly chosen loci.  {\it Genetics}.  {\bf 160}, 1179-1189.

\mn S. Sagitov (1999).  The general coalescent with asynchronous
mergers of ancestral lines.  {\it J. Appl. Probab.} {\bf 36}, 1116-1125.

\mn S. Sagitov (2003).  Convergence to the coalescent with simultaneous
multiple mergers.  {\it J. Appl. Probab.} {\bf 40}, 839-854.

\mn J. Schweinsberg (2000).  Coalescents with simultaneous multiple
collisions.  {\it Electron. J. Probab.} {\bf 5}, 1-50.

\mn J. Schweinsberg (2003).  Coalescent processes obtained from
supercritical Galton-Watson processes.  {\it Stochastic Process. Appl.}
{\bf 106}, 107-139.

\mn J. Schweinsberg and R. Durrett (2004).  Random partitions approximating
the coalescence of lineages during a selective sweep.  Preprint.
Available at http://front.math.ucdavis.edu/math.PR/ 0411069.

\mn K. Simonsen, G. A. Churchill, and C. F. Aquadro (1995).
Properties of statistical tests of neutrality for DNA polymorphism data.
{\it Genetics}, {\bf 141}, 413-429.

\mn W. Stephan, T. Wiehe, and M. W. Lenz (1992).  The effect of
strongly selected substitutions on neutral polymorphism: Analytical
results based on diffusion theory.  {\it Theor. Popul. Biol.} {\bf 41}, 237-254.

\mn F. Tajima (1989).  Statistical method for testing the neutral
mutation hypothesis by DNA polymorphism.  {\it Genetics}, {\bf 123}, 585-595.
\end{document}